\documentclass[reqno]{amsart}
\usepackage[T1]{fontenc}
\usepackage[utf8]{inputenc}
\usepackage[english]{babel}

\usepackage{blindtext}
\usepackage{amssymb}
\usepackage{hyperref}
\usepackage[dvipsnames]{xcolor}
\usepackage{tikz}
\usepackage{tikz-cd}
\usetikzlibrary{patterns}
\usepackage{xcolor}
\hypersetup{
    colorlinks,
    linkcolor={red!50!black},
    citecolor={blue!50!black},
    urlcolor={blue!80!black}
}
 
\usepackage{amsthm}
\numberwithin{equation}{section}
\usepackage{amsmath}
\usepackage{mathtools}
\usepackage{dirtytalk}
\usepackage{bookmark}

\DeclareFontFamily{U} {MnSymbolA}{}
\DeclareFontShape{U}{MnSymbolA}{m}{n}{
   <-6> MnSymbolA5
   <6-7> MnSymbolA6
   <7-8> MnSymbolA7
   <8-9> MnSymbolA8
   <9-10> MnSymbolA9
   <10-12> MnSymbolA10
   <12-> MnSymbolA12}{}
\DeclareFontShape{U}{MnSymbolA}{b}{n}{
   <-6> MnSymbolA-Bold5
   <6-7> MnSymbolA-Bold6
   <7-8> MnSymbolA-Bold7
   <8-9> MnSymbolA-Bold8
   <9-10> MnSymbolA-Bold9
   <10-12> MnSymbolA-Bold10
   <12-> MnSymbolA-Bold12}{}
\DeclareSymbolFont{MnSyA} {U} {MnSymbolA}{m}{n}
 
\DeclareMathSymbol{\top}{\mathord}{MnSyA}{219}
\DeclareMathSymbol{\bot}{\mathord}{MnSyA}{217}

\newtheorem{theorem}{Theorem}[section]
\newtheorem{prop}[theorem]{Proposition}
\newtheorem{corollary}[theorem]{Corollary}
\newtheorem{lemma}[theorem]{Lemma}

\DeclarePairedDelimiter{\ceil}{\lceil}{\rceil}

\theoremstyle{remark}
\newtheorem{remark}[theorem]{Remark}

\theoremstyle{question}
\newtheorem*{question}{Question}

\theoremstyle{remark}

\theoremstyle{definition}
\newtheorem{definition}[theorem]{Definition}
\newcommand{\norm}[1]{\left\lVert#1\right\rVert}

\DeclarePairedDelimiter\floor{\lfloor}{\rfloor}

\newcommand\blfootnote[1]{%
  \begingroup
  \renewcommand\thefootnote{}\footnote{#1}%
  \addtocounter{footnote}{-1}%
  \endgroup
}

\allowdisplaybreaks

\begin{document}

\title{The Lyapunov and uniform spectra}
\author{Jairo Bochi}
\address{Department of Mathematics, Penn State University}
\email{\href{mailto:bochi@psu.edu}{bochi@psu.edu}}
\author{Nicol\'as Mart\'inez Ramos}
\address{Department of Mathematics, Penn State University}
\email{\href{mailto:ndm5313@psu.edu}{ndm5313@psu.edu}}
\date{\today}

\begin{abstract}

The \textit{Lyapunov spectrum} of a linear cocycle (or vector bundle) is the set of Lyapunov vectors with respect to all ergodic probability measures  for the base topological dynamical system. This spectrum disregards trajectories that are atypical from the measure-theoretic viewpoint, and so it is in general insufficient to determine important properties of the cocycle, such as existence of dominated splittings. We resolve this deficiency by introducing a larger set which we call the \textit{uniform spectrum} of the cocycle. This definition is akin to Misiurewicz and Ziemian's rotation set. We prove several properties of the Lyapunov and uniform spectra, including relations between them. Motivated by a theorem of DeWitt and Gogolev, we consider Hölder-continuous cocycles over hyperbolic homeomorphisms and then provide bounds for the uniform spectrum in terms of the Lyapunov spectrum. We exhibit additional connections to the literature.  
\end{abstract}

\maketitle

\section{Introduction} 

\blfootnote{\copyright 2026 The authors. All rights reserved. It is hereby forbidden to upload the contents of this document, either in full or in part, to any generative AI system.}
Ergodic theory deals with measure-preserving dynamical systems. Nevertheless, this theory sometimes casts light on problems where no invariant measure is given a priori. In fact, ergodic theory has been used to answering questions from distant fields such as number theory. 

The classical theorem of Krylov and Bogolyubov states that every topological dynamical system (that is, a continuous action of $\mathbb{Z}^+$ or $\mathbb{R}^+$ on some compact metrizable space) admits at least one invariant probability measure. This foundational fact is to a great deal responsible for the applicability of ergodic theory to other areas.

Another basic fact, due to Oxtoby, is the following: if a topological dynamical system is uniquely ergodic (that is, admits a unique invariant probability measure),
then the Birkhoff averages of continuous functions converge uniformly. 
This is one of the simplest manifestations of an interplay between ergodic-theoretic and topological properties.

In this paper, we are interested in compositions of linear operators driven by some topological dynamical system. The asymptotic growth of those compositions, at least for typical orbits with respect to an invariant probability measure, is governed by the Lyapunov exponents, in the way described by the multiplicative ergodic theorem of Oseledets.

Walters \cite{Walters84} discovered that the theorem of Oxtoby mentioned above does \emph{not} extend to compositions of linear operators: unique ergodicity and continuity are not sufficient to ensure uniformity of  the limits defining the Lyapunov exponents.
There may exist trajectories along which the Lyapunov exponents are not well-defined, but these trajectories are not ``seen'' by the invariant measure.

In this paper, we consider vector bundle maps $\Phi$ covering a discrete-time topological dynamical system.
Thus we have a commutative diagram of the form 
\begin{equation}
\begin{tikzcd}
\mathbb{E} \arrow[r,"\Phi"] \arrow[d,swap,"\pi"] & \mathbb{E} \arrow[d,"\pi"] \\
X \arrow[r,swap,"T"] & X
\end{tikzcd}
\end{equation}
where $X$ is a compact metric space, $T$ is a continuous map, $\pi$ is the bundle projection, and $\Phi$ is a continuous map acting linearly on the fibers $\mathbb{E}_x \coloneq \pi^{-1}(x)$, which we assume are real vector spaces of finite dimension $d$.
In addition, we require that these fiber maps $\Phi_x \colon \mathbb{E}_x \to  \mathbb{E}_{Tx}$ are isomorphisms. Note that the iterates of $\Phi$ have fiber maps
\begin{equation}
	\Phi^n_x = \Phi_{T^{n-1}x} \circ \cdots \circ \Phi_x.
\end{equation}

This class of dynamical systems $\Phi$ includes tangent maps of local diffeomorphisms of compact manifolds, as well as linear cocycles over continuous maps.

Define the set 
\begin{equation}
\mathfrak{a}_{+}^{d}\coloneq\{{v}\in\mathbb{R}^{d}:{v}=(a_1,\dots,a_d) \text{ where }a_1\geq \dots\geq a_d\},
\end{equation}
which we call the (closed) \emph{positive Weyl chamber}. Given a bundle map $\Phi$ as above, we would like to determine all possible exponential growth rates of its iterates.
With this goal in mind, it is natural to consider the \emph{Lyapunov spectrum} $\mathcal{L}(\Phi)$, defined as the subset of $\mathfrak{a}^{d}_{+}$ consisting of all Lyapunov vectors of all ergodic invariant probability measures.
Thus, a vector $(\lambda_1,\dots,\lambda_d) \in \mathfrak{a}_{+}^d$ is an element of $\mathcal{L}(\Phi)$ if and only if there exists a $\Phi$-invariant ergodic probability measure $\mu$ on $X$ such that, for $\mu$-almost every $x$,
\begin{equation}
\lim_{n \to \infty} \frac{1}{n} \log \sigma_j(\Phi^n_x) = \lambda_j \, ,  \qquad j=1,\dots,d \, ,
\end{equation}
where $\sigma_j$ denotes the $j$\textsuperscript{th} largest singular value.

The Lyapunov spectrum depends only on the typical trajectories in the sense of ergodic theory, and thus may miss out relevant properties of $\Phi$. For example, a cocycle may have all Lyapunov exponents uniformly away from zero and yet fail to be uniformly hyperbolic: see \cite[Theorem~4.1]{Renato}. For examples using derivative cocycles, see \cite{CaoLuzzattoRios, Gogolev}. 

Therefore, we introduce a second subset $\mathcal{U}(\Phi)\subseteq\mathfrak{a}_{+}^d$, which we call \emph{uniform spectrum}, and consists of all vectors $(\xi_1,\dots,\xi_d)$ that can be obtained as limits of the form
\begin{equation}
\lim_{i \to \infty} \dfrac{1}{n_i}\big( \log \sigma_1(\Phi^{n_i}_{x_i}) , \dots ,  \log \sigma_d(\Phi^{n_i}_{x_i}) \big) \, ,
\end{equation}
where $(x_i)$ is a sequence of points in $X$, and $(n_i)$ is a sequence of positive integers tending to $+\infty$.
It is immediate that
\begin{equation}\label{e.trivial_inclusion}
	\mathcal{L}(\Phi) \subseteq \mathcal{U}(\Phi) \, . 
\end{equation}

While the Lyapunov spectrum $\mathcal{L}(\Phi)$, or variations of it, appears in the work of several authors \cite{DiazGelfertWangYang,Parkproximal,BreuillardSert}. Our definition of uniform spectrum $\mathcal{U}(\Phi)$ bears a similarity to the rotation set of a lift of a map of the torus homotopic to the identity, as defined by Misiurewicz and Ziemian \cite{MisiurewiczZiemian}. Other related objects include the Benoist cone \cite{Benoist} and the joint spectrum \cite{BreuillardSert}. 

Some uniform properties of the vector bundle map $\Phi$ can be read directly from the uniform spectrum $\mathcal{U}(\Phi)$. For example, $\Phi$ admits a dominated splitting if and only if $\mathcal{U}(\Phi)$ is disjoint from some chamber wall $\xi_{i}=\xi_{i+1}$ (see Proposition \ref{DominatedSplittingLemma}). Similarly, uniform hyperbolicity is detected by $\mathcal{U}(\Phi)$: see Proposition \ref{proposition uniform hyperbolicity}.

There are \emph{some} uniform properties of $\Phi$ that can be read from the Lyapunov spectrum $\mathcal{L}(\Phi)$. 
In fact, the Lyapunov spectrum imposes restrictions on the uniform spectrum: more precisely, 
\begin{equation}\label{e.reverse_inclusion}
	\mathcal{U}(\Phi) \subseteq \mathrm{Maj}(\mathcal{L}(\Phi)) \, ,
\end{equation}
where the set $\mathrm{Maj}(S)$ denotes the \textit{majorization hull} of a subset $S\subseteq \mathfrak{a}_{+}^{d}$, defined as follows
\begin{equation}
    \mathrm{Maj}(S)\coloneq \left\{\xi\in\mathfrak{a}_{d}^{+} \ :\forall b\in\mathfrak{a}^{d}_{+}, \ \langle b,\xi\rangle\leq \sup_{\eta\in S} \langle b,\eta\rangle  \right\},
\end{equation}
where  $\langle\cdot,\cdot\rangle$ is the standard inner product on $\mathbb{R}^{d}$. The majorization hull is discussed in detail in Section \ref{Section 3}. Some examples are presented in Figures \ref{Majorization hull two dimensional example} and \ref{Majorization hull SL(3) example}. Inclusion (\ref{e.reverse_inclusion}) is Corollary \ref{First bound of uniform spectrum}. 

It is unclear whether there are any other general relations between out two spectra other than the inclusions \eqref{e.trivial_inclusion} and \eqref{e.reverse_inclusion}. On the other hand, we do have extra information for specific classes of bundle maps $\Phi$.

Our next result relies on the following definition. If $\gamma\geq 1$ is a parameter, the \textit{$\gamma$-aura} of a bounded set $S\subseteq \mathfrak{a}^{d}_{+}$ is 
\begin{align}
\begin{split}
    \mathrm{Aura}_{\gamma}(S)\coloneq &\left\{\xi \in  \mathfrak{a}_{d}^{+} \ :\right. \\
    &\left. \forall b\in\mathfrak{a}^{d}_{+}, \ (1-\gamma)\sup_{\eta\in S}\langle b,\eta\rangle+\gamma\inf_{\eta\in S} \langle b,\eta\rangle \leq\langle b,\xi\rangle\leq \sup_{\eta\in S} \langle b,\eta\rangle  \right\}.
\end{split}
\end{align}
It is clear that Aura$_{\gamma}(S)\subseteq \mathrm{Maj}(S)$. See Figures \ref{SwellingExample} and \ref{Eye} for examples. Theorem \ref{MainTheorem}, which is a major result of this paper, states that, if $T$ is a topologically transitive hyperbolic homeomorphism and $\Phi$ is a Hölder bundle map, then 
\begin{equation}
    \mathcal{U}(\Phi)\subseteq \mathrm{Aura}_{\gamma_0}(\mathcal{L}(\Phi)), 
\end{equation}
where $\gamma_{0}\geq 1$ is an explicit parameter depending only on the hyperbolicity rates of $T$, the Hölder exponent of $\Phi$, and on a certain measurement of the non-conformality of $\Phi$. In case $\Phi$ satisfies the well known \textit{fiber-bunching} condition, the parameter $\gamma_{0}$ equals $1$.

Let us demonstrate the relevance of Theorem \ref{MainTheorem} by deriving some results from the literature as corollaries. Velozo-Ruiz \cite[Theorem~1]{Renato} showed that for a fiber-bunched 2-dimensional linear cocycle over a hyperbolic base admits a dominated splitting if there is a uniform gap between its Lyapunov exponents. The latter condition means that the Lyapunov spectrum $\mathcal{L}(\Phi)$ is away from the chamber wall $\xi_{1}=\xi_{2}$. In dimension $d=2$, the $1$-aura of any bounded subset of the half-space $\mathfrak{a}_{+}^{2}$ equals the closed convex hull $\overline{\mathrm{conv}}(S)$: see Proposition \ref{Aura in dimension 2}. Therefore, Theorem \ref{MainTheorem} states that in this situation 
\begin{equation}
    \mathcal{U}(\Phi)\subseteq \overline{\mathrm{conv}}(\mathcal{L}(\Phi)).
\end{equation}
Now, Velozo-Ruiz's theorem follows from the aforementioned characterization of dominated splittings in terms of the uniform spectrum. 

As a second application of Theorem \ref{MainTheorem} we will reobtain a result of DeWitt and Gogolev \cite[Theorem~1.3]{DewittGogolev}, which actually motivated our investigation. Let $\Phi$ be a Hölder bundle map over a hyperbolic base. Suppose $\mathcal{L}(\Phi)$ is contained in a sufficiently small ball around a point $\xi\in\mathfrak{a}^{d}_{+}$ such that $\xi_{i}>\xi_{i+1}$ for some $i\in \{1,\dots, d-1\}$. Then the theorem of DeWitt and Gogolev states that $\Phi$ admits a dominated splitting of index $i$. This statement is a corollary of Theorem \ref{MainTheorem} since the auras of small sets are small. Therefore, our result is a qualitative and quantitative refinement of theirs.  

\medskip

\noindent\pdfbookmark[0]{Structure}{Structure}\textsc{Organization of the paper.} In Section \ref{Section 2}, we formally define the two spectra $\mathcal{L}(\Phi)$ and $\mathcal{U}(\Phi)$ and state several general properties. In Section \ref{Section 3}, we recall the classical majorization partial order and the notion of Schur-convexity,  which we use to study majorization hulls and auras. Section \ref{Section 4} states and proves Theorem \ref{MainTheorem}, whose main applications are detailed in Section \ref{Corollaries}. In Section \ref{Section questions}, we review some results of Park \cite{Parkquasimultiplicativity,Parkproximal} in terms of our two spectra and we state a couple of questions. The proofs of some auxiliary statements in this paper use arguments dissimilar to the main ideas presented; these proofs were moved to Appendices \ref{properties Lyapunov and uniform spectra} and \ref{Appendix B}.



\bigskip

\noindent\pdfbookmark[0]{AI statement}{AIstatement}\textsc{Disclosure notice.} All mathematical definitions, statements and proofs in this paper were formulated by the authors without assistance of generative artificial intelligence. No Gen AI was used in the literature search, nor in the preparation of the figures. The LLMs embedded in Overleaf and VS Code were employed in correcting grammar and minor language polishing. 

\bookmarksetup{startatroot}

\section{The two spectra} \label{Section 2}

\subsection{The Lyapunov spectrum}

Throughout this paper we will assume that $(X,d)$ is a compact metric space and $T:X\to X$ is a continuous map. Furthermore, we will assume that $\mathbb{E}$ is a $d$-dimensional $C^0$ vector bundle over $X$ with projection $\pi:\mathbb{E}\to X$. More precisely, there is a finite open cover $\{U_i\}$ of $X$  and local homeomorphisms $\psi_{i}:U_i\times \mathbb{R}^{d}\to \pi^{-1}(U_i)$ such that for all $(x,u)\in U_i\times \mathbb{R}^{d}$ we obtain $\pi(\psi(x,u))=x$. Additionally,  for $x\in U_i\cap U_j$ the map $[\psi_{j}(x,\cdot)]^{-1}\circ \psi_{i}(x,\cdot)$ is linear. We endow each fiber $\mathbb{E}_{x}=\pi^{-1}\{x\}$ with the structure of a $d-$dimensional such that 
\begin{equation} \label{definition of h}
  h_{i}(x):=\psi_i(x,\cdot):\mathbb{R}^{d}\to\mathbb{E}_x  
\end{equation}
is an isomorphism.

We call $\Phi:\mathbb{E}\to\mathbb{E}$ a \textit{bundle map} of $\mathbb{E}$ if $\Phi$ covers the map $T$, that is, $\pi\circ\Phi=T\circ\pi$  and if for every $x\in X$ the restricted map $\Phi_{x}:\mathbb{E}_{x}\to \mathbb{E}_{Tx}$ is linear. For this paper, we will also assume $\Phi$ is non-singular, that is, every $\Phi_{x}$ is an isomorphism. 

 We will assume that $\mathbb{E}$ is endowed with a $C^0$ \textit{Riemannian metric} $g$, that is, a continuous family of inner products $g_{x}=\langle\cdot,\cdot\rangle_{x}$ where $x\in X$.


Given two points $x,y\in X$, denote the set of invertible linear maps from $\mathbb{E}_{x}$ to $\mathbb{E}_{y}$ as $\text{GL}(\mathbb{E}_{x},\mathbb{E}_{y} )$. Let $A$ be in $\text{GL}(\mathbb{E}_{x},\mathbb{E}_{y})$. The operator norm is defined as
\begin{equation}
 \norm{A}_{y\xleftarrow{}x}=\sup_{\substack{v\in\mathbb{E}_{x}\\v\not=0}}\dfrac{\norm{Av}_y}{\norm{v}_x}.
 \end{equation}
 We omit the subscripts when they are clear from the context.
 
There exists $A^{*}\in \text{GL}(\mathbb{E}_{y},\mathbb{E}_{x})$ such that for all $v\in \mathbb{E}_{x}$, $w\in \mathbb{E}_{y}$
\begin{equation} \label{adjoint matrix}
  \langle w, Av\rangle_{x}=\langle A^{*}w, v\rangle_{y}.  
\end{equation}
Define the \textit{singular values} $\sigma_{1}(A)\geq\dots\geq\sigma_{d}(A)$ of $A$ as the eigenvalues of the positive self-adjoint operator $A^{*}A\in \text{GL}(\mathbb{E}_{x})$. The first singular value $\sigma_{1}(A)$ equals the norm $\norm{A}$. The last singular value $\sigma_{d}(A)$ is called the \textit{conorm} of $A$ and denoted $\mathfrak{m}(A)$. Note that  $\mathfrak{m}(A)=\norm{A^{-1}}^{-1}$. More generally, for $i\in\{1,\dots,d\}$

\begin{equation}
    \sigma_{i}(A)=[\sigma_{d+1-i}(A^{-1})]^{-1}.
\end{equation}

Define the \textit{Cartan projection} as $\kappa^{g}:\text{GL}(\mathbb{E}_{x},\mathbb{E}_{y} )\to\mathfrak{a}_{+}^{d}$ as 
\begin{equation}
    \kappa^{g}(A)\coloneq(\log \sigma_1(A),\dots, \log \sigma_d(A)).
\end{equation}
The superscript $g$ is a reminder that the Cartan projection depends on the choice of Riemannian metric. 

The \textit{Jordan projection} is defined as $\zeta:\text{GL}(\mathbb{E}_{x},\mathbb{E}_x)\to \mathfrak{a}_{+}^{d}$ as
\begin{equation}
 \zeta(A)\coloneq(\log|\chi_1(A)|,\dots, \log|\chi_d(A)|),
\end{equation}
where $|\chi_1|,\dots, |\chi_d|$ are the absolute values of the eigenvalues $A$ in decreasing order. The two projections are related by the formula 
\begin{equation}\label{Gelfand}
    \zeta(A)=\lim_{n\to\infty} \dfrac{1}{n}\kappa^{g}(A^{n}),
\end{equation}
which is a generalization of Gelfand's formula $|\chi_{1}(A)|=\lim_{n\to\infty}\norm{A^{n}}^{1/n}$.

Denote by $\mathcal{M}_{\mathrm{erg}}(X)$ the set of Borel probability measures which are ergodic with respect to $T$. Let $\mu\in \mathcal{M}_{\mathrm{erg}}(X)$. By the Furstenberg-Kesten theorem, the sequence of functions $x\mapsto \frac{1}{n}\kappa^{g}(\Phi^{n}x)$ converges $\mu$-almost everywhere and in $\mathrm{L}^{1}(\mu)$-norm to a constant vector in $\mathfrak{a}_{+}^{d}.$


\begin{definition}
    The \textit{Lyapunov spectrum} of the bundle map $\Phi:\mathbb{E}\to \mathbb{E}$ is the set
    \begin{align}
    \begin{split}
    \mathcal{L}(\Phi)\coloneq &\left\{\xi \in \mathfrak{a}_{+}^{d}:\exists\nu\in \mathcal{M}_{\mathrm{erg}}(X) \text{ such that }\right.\\
    &\left. \hspace{3.5cm}\xi=\lim_{n\to\infty}\dfrac{1}{n}\int\kappa^{g}(\Phi^{n}_{x})d\nu \ \right\}.
    \end{split}
    \end{align}
\end{definition}
This set is independent of the Riemannian metric $g$. Since $X$ is a compact space and $T$ is a continuous transformation, and $\Phi$ is nonsingular, $\mathcal{L}(\Phi)$ is nonempty and bounded. Being the image of the Borel space $\mathcal{M}_{\mathrm{erg}}(X)$ under a Borel map, the Lyapunov spectrum $\mathcal{L}(\Phi)$ is an analytic (that is, Souslin) set.

\begin{prop}\label{PeriodicPointsLyapunovSpectrum}
Each periodic orbit of $T$ gives rise to an element of $\mathcal{L}(\Phi)$ as follows:
    \begin{equation}
    \left\{\frac{1}{n}\zeta(\Phi^{n}_{p}):\ n\in\mathbb{N},\ T^{n}p=p\right\}\subseteq \mathcal{L}(\Phi).
    \end{equation}
\end{prop}
\begin{proof}
This is  consequence of formula (\ref{Gelfand}) and the fact that every periodic orbit is the support of an ergodic measure.  
\end{proof}
\subsection{The uniform spectrum} 

\begin{definition}
    The \textit{uniform spectrum} of $\Phi$ is the set
    \begin{align}
    \begin{split}
      \mathcal{U}(\Phi)\coloneq &\left\{\xi \in\mathfrak{a}_{+}^{d}:\exists \{n_j\}\subseteq\mathbb{N} \text{ where }n_{j}\to\infty \text{, } \{x_j\}\subseteq X \right.\\
      &\left. \hspace{3.5cm}\text{ such that }\xi =\lim_{j\to\infty}\dfrac{1}{n_j}\kappa^{g}(\Phi_{x_j}^{n_j})\right\}.
    \end{split}
    \end{align}
\end{definition}

Note that $v\in \mathcal{U}(\Phi)$ if and only if for any neighborhood $\mathcal{O}\subseteq \mathfrak{a}_{+}^{d}$ of $v$ there exists a sequence $\{x_n\}\subseteq X$ such that 
$n^{-1}\kappa^{g}(\Phi_{x_n}^{n})\in \mathcal{O}$ infinitely often. Hence, if for all $n\in\mathbb{N}$ we write
\begin{equation}
\kappa^{g}(\Phi^{n}):=\{\kappa^{g}(\Phi^{n}_{x}):x\in X\},
\end{equation}
then we can characterize the uniform spectrum as
\begin{equation}\label{KuratowskiLimsup}
\mathcal{U}(\Phi)=\underset{{n\to\infty}}{\text{Ls}}\dfrac{\kappa^{g}(\Phi^n)}{n},
\end{equation}
where $\text{Ls}$ is the Kuratowski limit superior \cite[page~337]{Kuratowski}.

Like $\mathcal{L}(\Phi)$, the set $\mathcal{U}(\Phi)$ is independent of the Riemannian metric $g$. Moreover, from its characterization  as a Kuratowski limit,  $\mathcal{U}(\Phi)$ is a compact set. Furthermore, it is clear that 
\begin{equation}
  \mathcal{L}(\Phi)\subseteq \mathcal{U}(\Phi).  
\end{equation}
Therefore, the set $\mathcal{U}(\Phi)$ is nonempty.

On both the Lyapunov and uniform spectra, we can find the following relation between the spectra of $\Phi$ and $\Phi^{k}$ for any $k\in\mathbb{N}$
\begin{prop} \label{change of uniform spectrum due to iterations}
    For every $k\in \mathbb{N}$, 
    \begin{equation} \label{Lyapunov spectrum iterations}
        \mathcal{L}(\Phi^{k})=k\ \mathcal{L}(\Phi);
    \end{equation}
    \begin{equation}\label{Uniform spectrum iterations}
        \mathcal{U}(\Phi^{k})=k\ \mathcal{U}(\Phi) .
    \end{equation}
    Furthermore, if $T$ is invertible, then 
    \begin{align}
        \mathcal{L}(\Phi^{-1})=\iota(\mathcal{L}(\Phi)); \label{Inverse Lyapunov Spectrum}\\
        \mathcal{U}(\Phi^{-1})=\iota(\mathcal{U}(\Phi)), \label{inverse uniform spectrum}
    \end{align}
    where $\iota:\mathfrak{a}_{+}^{d}\to\mathfrak{a}_{+}^{d}$, $\iota(\xi_{1},\dots, \xi_{d})\coloneq(-\xi_{d},\dots, -\xi_{1})$ is the opposition involution. 
\end{prop}

Similar properties of rotation sets are given in \cite[Proposition~2.1]{MisiurewiczZiemian}. We present the proof of this proposition to Appendix \ref{properties Lyapunov and uniform spectra}.

\subsection{Relation between the uniform spectrum and domination} \label{UniformSpectrum}
    In this subsection we will show how the existence of dominated splittings can be determined from the uniform spectrum  
    \begin{definition}
        A continuous $\Phi$-invariant splitting $\mathbb{E}=\mathbb{E}^{1}\oplus \mathbb{E}^{2}$ is \textit{dominated} for $\Phi$ if it is invariant under $\Phi$ and there are constants $C>0$ and $\tau\in(0,1)$ such that for every $x\in X$ and $n\in\mathbb{N}$
        \begin{equation} \label{definicion de m}
        \dfrac{\norm{\Phi^{n}_{x}|\mathbb{E}^{2}_{x}}}{\mathfrak{m}(\Phi^{n}_{x}|\mathbb{E}^{1}_{x})}\leq C\tau^{n}.
        \end{equation}
        The dimension of $\mathbb{E}^{1}$ is called \textit{the index of domination} of the splitting. \\
        It can be shown (see \cite[Appendix~B.1]{BonattiDiazViana}) that for any $\Phi$ there exists a unique continuous $\Phi$-invariant splitting 
        \begin{equation}
          \mathbb{E}=\mathbb{E}^{1}\oplus\dots\oplus\mathbb{E}^{m}\ \text{possibly with }m=1,  
        \end{equation}
        called the \textit{finest dominated splitting}, such that all dominated splittings into two invariant subbundles are of the form $(\mathbb{E}^{1}\oplus\dots\oplus\mathbb{E}^{k})\oplus (\mathbb{E}^{k+1}\oplus\dots\oplus\mathbb{E}^{m}).$
    \end{definition}

    \begin{prop} \label{DominatedSplittingLemma}
        Assume that $T:X\to X$ is a continuous map and $\Phi:\mathbb{E}\to\mathbb{E}$ is a bundle map covering $T$. Then there is a dominated splitting of index $i\in \{1,\dots, d-1\}$ if and only if
        \begin{equation}
         \min_{\xi\in \mathcal{U}(\Phi)}(\xi_i-\xi_{i+1})>0,   
        \end{equation}
        that is, the set $\mathcal{U}(\Phi)$ does not intersect  the i$^{th}$ wall of the positive Weyl chamber. 
    \end{prop}
    \begin{proof}
        According to \cite[Theorem~A]{BochiGourmelon}, $\Phi$ admits a dominated splitting of index $i$ if and only if there exist $n_{0}\in\mathbb{N}$ and $\epsilon>0$ such that for all $x\in X$ and $n>n_0$, 
    \begin{equation} 
        \dfrac{\sigma_{i}(\Phi_{x}^n)}{\sigma_{i+1}(\Phi_{x}^n)}>e^{\epsilon n}.
    \end{equation}
    Equivalently, the sets $\kappa^{g}(\Phi^{n})/n$ are, for sufficiently large $n$, uniformly away from the wall $\xi_{i}=\xi_{i+1}$. This means that $\mathcal{U}(\Phi)$ does not intersect the wall.
    \end{proof}
   By inspecting the intersection between the boundary of $\mathfrak{a}_{+}^{d}$ and the uniform spectrum it is possible to obtain all the indexes of domination, or equivalently, the dimensions of the subbundles of the finest dominated splitting. For example, $\mathcal{U}(\Phi)\subseteq \mathrm{int}(\mathfrak{a}_{+}^{d})$ if and only if the finest dominated splitting is formed by $d$ one-dimensional subbundles. 

    Beyond domination, the uniform spectrum also yields information about the hyperbolicity of the system. 

    \begin{prop} \label{proposition uniform hyperbolicity}
        Assume that $T:X\to X$ is a continuous map and $\Phi:\mathbb{E}\to\mathbb{E}$ is a bundle map. 
        \begin{enumerate}
            \item $\Phi$ is uniformly contracting if and only if $\max_{\xi\in\mathcal{U}(\Phi)}\xi_{1}<0$. 
            \item $\Phi$ is uniformly expanding if and only if $\min_{\xi\in\mathcal{U}(\Phi)}\xi_{d}>0$.
            \item $\Phi$ has a uniformly hyperbolic splitting into two nontrivial bundles if and only if there exist $i\in\{1,\dots, d-1\}$ such that
            \begin{equation} \label{UniformHyperbolicity}
                \min_{\xi\in\mathcal{U}(\Phi)}\xi_{i}>0>\max_{\xi\in\mathcal{U}(\Phi)}\xi_{i+1}.
            \end{equation}
        \end{enumerate}
    \end{prop}

    \begin{proof} We start with the proof of $(1)$. Suppose that $\Phi$ is uniformly contracting. Then there exist $0 <\alpha< 1$ and $C>0$ such that for all $x\in X$ and all $n\in \mathbb{N}$
            \begin{equation}
                \norm{\Phi^{n}_{x}}\leq C\alpha^{n}.
            \end{equation}
            Thus, taking any sequences $\{x_k\}\subseteq X$ and $\{n_k\}\subseteq \mathbb{N}$ such that $n_k\to \infty$, we obtain 
            \begin{equation}
                \limsup_{k\to\infty}\dfrac{1}{n_k}\log \norm{\Phi^{n_k}_{x_k}}\leq \log \alpha. 
            \end{equation}
            Hence, for all $\xi\in \mathcal{U}(\Phi)$, we obtain that $\xi_{1}<0.$
            
            Now, suppose that $\Phi$ is not uniformly contracting. Then, for every $0<\alpha<1$, for every $C>0$ and for every $n\in\mathbb{N}$ it is possible to find a sequence $\{x_{n}\}\subseteq X$ such that 
            \begin{equation}
                \norm{\Phi^{n}_{x_n}}>C\alpha^{n}. 
            \end{equation}
            Passing to a converging subsequence given by the points $\{x_{k}\}\subseteq\{x_{n}\}$ such that $k\to\infty$, we obtain the inequality
            \begin{equation}
                \dfrac{1}{k}\log\norm{\Phi^{k}_{x_k}}>\dfrac{1}{k}\log C+\log\alpha. 
            \end{equation}
            Taking $k\to\infty$, there exist some $\eta=(\eta_1,\dots, \eta_{d})\in \mathcal{U}(\Phi)$ such that 
            \begin{equation}
                \eta_{1}=\lim_{k\to\infty}\dfrac{1}{k}\log\norm{\Phi^{k}_{x_k}}>\log\alpha.
            \end{equation}
            Since $\alpha$ is any number in the interval $(0,1)$, statement (1) follows. 
            
            We will now prove $(2)$.  Note that $\Phi$ is uniformly expanding if  there exist $0 <\alpha< 1$ and $C>0$ such that 
            \begin{equation}
                \norm{(\Phi^{n}_{x})^{-1}}\leq C\alpha^{n}.
            \end{equation}
            Hence, the argument used for (1) yields the result. 
            
            Finally, we prove $(3)$. Suppose $\Phi$ has a uniformly hyperbolic splitting into nontrivial bundles $\mathbb{E}=\mathbb{F}^{1}\oplus \mathbb{F}^{2}$. That is, $\Phi|_{\mathbb{F}^{1}}$ is uniformly expanding and $\Phi|_{\mathbb{F}^{2}}$ is uniformly contracting.  Taking $i=\dim \mathbb{F}^{1}$, from (1) and (2) we obtain that $\xi_{i}>0>\xi_{i+1}$ for all $\xi\in \mathcal{U}(\Phi)$. 

            On the other hand, if there is $i\in \{1,\dots, d-1\}$ such that (\ref{UniformHyperbolicity}) is obtained, by Proposition \ref{DominatedSplittingLemma}, there is an invariant splitting $\mathbb{E}=\mathbb{F}^{1}\oplus \mathbb{F}^{2}$ of index $i$ where $\mathbb{F}^{1}$ is uniformly expanding by (2) and $\mathbb{F}^{2}$ is uniformly contracting by (1). Statement (3) follows.  
    \end{proof}

\begin{remark}
    The Sacker-Sell spectrum is defined as 
    \begin{equation}
        \mathrm{SS}(\Phi)\coloneq\{\lambda\in\mathbb{R}: e^{-\lambda}\Phi \text{ is not uniformly hyperbolic}\}.
        \end{equation}
    It follows from Proposition \ref{proposition uniform hyperbolicity} that the set $\mathrm{SS}(\Phi)$ is determined by $\mathcal{U}(\Phi)$ and is actually a finite union of compact intervals, recovering a result from \cite{SackerSell}. We omit the details. 
\end{remark}

\subsection{A connectedness criterion}

We will show that under a mild condition on the dynamics, the uniform spectrum is a connected set.

Following \cite[p.~3]{LaSalle}, we say that a continuous map $T \colon X \to X$ acting on a compact metric space $X$ is \emph{invariantly connected} if $X$ is not the union of two disjoint nonempty closed forward-invariant sets. 
This property clearly holds if $X$ is connected or if $T$ is transitive; it also holds if $T$ is uniquely ergodic, for instance.

\begin{prop}\label{p.criterion_U_connected}
If $T$ is invariantly connected, then for any  bundle map 
$\Phi$ covering $T$, the uniform spectrum $\mathcal{U}(\Phi)$ is connected.
\end{prop}

The proof of this proposition will be presented in Appendix \ref{properties Lyapunov and uniform spectra}. Compare with \cite[Theorem~2.3]{MisiurewiczZiemian}. 

Let us observe that the corresponding statement for $\overline{\mathcal{L}(\Phi)}$ is not true. Indeed, letting $T$ be the doubling map on the circle, the paper \cite{Bochi26} constructs a non-singular bundle map $\Phi$ covering $T$ such that $\mathcal{L}(\Phi)$ has an isolated point. We remark that the set $\mathcal{M}_\mathrm{erg}(T)$ of ergodic measures for the doubling map is path-connected; this fact follows by a minor adaptation of \cite[Theorem~B]{Sigmund77}. Other than unique ergodicity, we do not know any condition on the dynamics that will guarantee that $\mathcal{L}(\Phi)$ will be connected for all non-singular bundle maps $\Phi$ covering $T$.


\subsection{Spectra of $\mathrm{SL}(2,\mathbb{R})$-cocycles}

To illustrate the definitions we have given so far, let us discuss the case where the fiber dimension of $\mathbb{E}$ is $d=2$ and $|\det \Phi_{x}|=1$ for all $x\in X$ (where $|\det|$ denotes the product of the singular values). This class of bundle maps $\Phi$ includes all continuous $\mathrm{SL}(2,\mathbb{R})$-cocycles. Then both spectra $\mathcal{L}(\Phi)$ and $\mathcal{U}(\Phi)$ are contained in the half-line $\{(t,-t):t\geq 0\}\subseteq \mathfrak{a}^{d}_{+}$, which we identify with the interval $[0,\infty)$ by projecting on the first coordinate. Basically without loss of generality, we assume that $T$ is invariantly connected. Then, by Proposition \ref{p.criterion_U_connected}, the uniform spectrum $\mathcal{U}(\Phi)$ is connected (and compact), hence equal to some interval $[\alpha^{\bot},\alpha^{\top}]$. By Proposition \ref{proposition uniform hyperbolicity}, $\alpha^{\bot}>0$ if and only if $\Phi$ is uniformly hyperbolic. As for the Lyapunov spectrum $\mathcal{L}(\Phi)$, we know it is contained in the interval $[\alpha^{\bot},\alpha^{\top}]$. By the forthcoming Proposition \ref{proposition about top}, we actually have $\alpha^{\top}\in \mathcal{L}(\Phi)$. If $\Phi$ is uniformly hyperbolic, restricting $\Phi$ to the unstable bundle it can be shown that $\alpha^{\bot}\in \mathcal{L}(\Phi).$

Suppose $T$ is uniquely ergodic. Then $\mathcal{L}(\Phi)=\{\alpha^{^{\top}}\}$. Furthermore, $\alpha^{\bot}$ equals either $0$ or $\alpha^{\top}$ (as a consequence of the hyperbolic subbundles if $\alpha^{\bot}>0$). If $\alpha^{\bot}=0$, then for every $x$ in some residual subset of $X$, the set of accumulation points of the sequence $\{\frac{1}{n}\kappa^{g}(\Phi^{n}_{x})\}_{n\in\mathbb{N}}$ is exactly $\mathcal{U}(\Phi)$; this is the content of \cite[Corollary~4.13]{BPS26}. 

Next, assume that $T$ is a subshift of finite type. The example constructed by Velozo-Ruiz \cite[Theorem~5]{Renato} shows that it is possible to have both $\alpha^{\bot}=0$ and $\mathcal{L}(\Phi)\subseteq [\epsilon, \alpha^{\top}]$ for some $\epsilon>0$ (even if $\Phi$ is a Hölder cocycle). The paper \cite{Monochromatic} constructs continuous cocycles $\Phi$ such that $\alpha^{\bot}=0<\alpha^{\top}$ and $\mathcal{L}(\Phi)=\{\alpha^{\top}\}$; however, these cocycles cannot be Hölder by \cite[Theorem~1.3]{DewittGogolev}. 

\section{Majorization} \label{Section 3}
\subsection{The majorization partial order}
We define the following partial order on $\mathbb{R}^{d}$. Let $\mathrm{S}_{d}$ be the groups of permutations of $\{1,\dots, d\}$. 
\begin{definition} \label{majorization}
    Consider two vectors $\xi, \eta\in \mathbb{R}^{d}$. Let $\tau,\sigma\in \mathrm{S}_d$ be permutations of $d$ elements such that $(\xi_{\tau(1)},\dots,\xi_{\tau(d)}), (\eta_{\sigma(1)},\dots,\eta_{\sigma(d)})\in \mathfrak{a}_{+}^{d} $. We will say ${\eta}$ \textit{majorizes} ${\xi}$ and we will write ${\xi} \trianglelefteq {\eta}$ if for all $k\in\{1,\dots, d-1\}$
    \begin{align}
    \sum_{i=1}^{k}\xi_{\tau(i)}\leq \sum_{i=1}^{k}\eta_{\sigma(i)} \ \ \text{ and }\ \ \sum_{i=1}^{d}\xi_{i}= \sum_{i=1}^{d}\eta_{i}.
    \end{align}
\end{definition}

Weyl's majorization theorem \cite[Theorem~9.E.1]{majorization}, gives a relation between the Cartan and Jordan projections: for any $A\in\text{GL}(\mathbb{E}_{x},\mathbb{E}_{x})$,
\begin{equation} \label{EigenvaluesVSSingularValues}
        \zeta(A)\trianglelefteq \kappa^{g}(A). 
\end{equation}
We also have the following subadditivity property \cite[Theorem~9.H.1]{majorization}: Given  $A\in \text{GL}(\mathbb{E}_{x},\mathbb{E}_{y})$ and $B\in \text{GL}(\mathbb{E}_{y},\mathbb{E}_{z})$ for $x,y,z\in X$, then 
\begin{equation} \label{TriangleInequalityCartanProjection}
\kappa^{g}(BA)\trianglelefteq\kappa^{g}(B)+\kappa^{g}(A).
\end{equation}

The following is a useful geometric characterization of majorization \cite[Corollary~2.B.3]{majorization}.

\begin{prop}  \label{MajorizationByConvexHull}
    The set $\{\xi\in\mathbb{R}^{d}: \xi\trianglelefteq \eta\}$ is the permutohedron obtained as the convex hull of the set $\{(\eta_{\sigma(1)},\dots,\eta_{\sigma(d)}):\sigma\in \mathrm{S}_d\}$. 
\end{prop}

\begin{remark}  \label{TrivialMajorization}
   Vectors of the form $(t,\dots, t)$ are minimal elements of the partial order. Every $\xi\in\mathbb{R}^{d}$ majorizes a unique vector of this form, namely with $t=\frac{1}{d}\sum_{i=1}^{d}\xi_{i}$.  
\end{remark}
 
For more information on majorization, we direct the reader to the book \cite{majorization}. 

    \subsection{The majorization hull}
    
    \begin{definition}
        A set $C\subseteq \mathfrak{a}_{+}^{d}$ is \textit{LCM-closed} if it satisfies the following three properties:
        \begin{enumerate}
        \item It is closed. 
        \item It is convex.
        \item If $\xi\trianglelefteq \eta$ for some $\eta\in C$, then $\xi\in C.$
    \end{enumerate}
    The prefix LCM represents the three ways these sets are closed: by taking limits, by taking convex conbinations and by majorization. 
    \end{definition}

    The following proposition gives a characterization of LCM-closed sets in terms of the action of the permutation group $\mathrm{S}_d$ on $\mathbb{R}^{d}.$

    \begin{prop} \label{MconvexProposition}
    Let $C\subseteq \mathfrak{a}_{+}^{d}$ be a closed set. Then $C$ is LCM-closed if and only if the set $\mathrm{S}_d\cdot C\subseteq \mathbb{R}^{d}$ is convex. 
\end{prop}
\begin{proof}
    First, let us show that if $C\subseteq \mathfrak{a}_{+}^{d}$ is LCM-closed then $\mathrm{S}_d\cdot C$ is convex. Consider $v,w\in \mathrm{S}_d\cdot C$ and $u=tv+(1-t)w$ for $t\in [0,1]$. By applying an appropriate permutation, without loss of generality we may assume that $u\in \mathfrak{a}_{+}^{d}.$ We have to show that $u\in \mathrm{S}_d\cdot C.$
    
    Notice that there are  $\tau,\sigma\in \mathrm{S}_d$ and $\xi,\omega\in C$ such that $v=\tau(\xi)=(\xi_{\tau(1)},\dots, \xi_{\tau(d)})$ and $w=\sigma(\omega)=(\omega_{\sigma(1)},\dots, \omega_{\sigma(d)})$. Consider $\eta=t\xi+(1-t)\omega\in C$. Let us show that $u\trianglelefteq \eta$. Since $\xi,\omega\in\mathfrak{a}_{+}^{d}$, we see that for $1\leq k\leq d$
    \begin{align}
    \begin{split}
        \sum_{i=1}^{k}\xi_{\tau(i)}&\leq \sum_{i=1}^{k}\xi_{i},\\
        \sum_{i=1}^{k}\omega_{\sigma(i)}&\leq \sum_{i=1}^{k}\omega_{i},
    \end{split}
    \end{align}
    with equality in the case that $k=d$. Hence, we obtain that
    \begin{align}
    \begin{split}
        \sum_{i=1}^{k}u_{i}&=t\sum_{i=1}^{k}v_{i}+(1-t)\sum_{i=1}^{k}w_{i}\\
        &=t\sum_{i=1}^{k}\xi_{\tau(i)}+(1-t)\sum_{i=1}^{k}\omega_{\sigma(i)}\\
        &\leq t\sum_{i=1}^{k}\xi_{i}+(1-t)\sum_{i=1}^{k}\omega_{i}\\
        &=\sum_{i=1}^{k} \eta_{i},
    \end{split}
    \end{align}
    with equality in the case that $k=d$. We conclude that $u\trianglelefteq \eta$, which implies that $u\in C$ since $C$ is LCM-closed.

    Now let us show that $C\subseteq \mathfrak{a}_{+}^{d}$ is LCM-closed if $C$ is closed and $\mathrm{S}_d\cdot C$ is convex. We see that $C=\mathrm{S}_d\cdot C\cap \mathfrak{a}_{+}^{d}.$ Since $\mathfrak{a}_{+}^{d}$ is convex, we can see that $C$ is convex. Now consider $\eta\in C$ and $\xi\in\mathfrak{a}_{+}^{d}$ such that $\xi\trianglelefteq \eta$. Hence, by Proposition \ref{MajorizationByConvexHull} we know that $\xi\in\mathrm{conv}\{(\eta_{\sigma(1)},\dots,\eta_{\sigma(d)}):\sigma\in \mathrm{S}_d\}\subseteq \mathrm{S}_d\cdot C$. We conclude that $\xi\in C$. 
\end{proof}

\begin{definition}
    Given any subset $S\subseteq \mathfrak{a}_{+}^{d}$, the \textit{majorization hull} of $S$, denoted $\mathrm{Maj}(S)$, is the smallest LCM-closed subset of $\mathfrak{a}^{d}_{+}$ that contains $S$.
\end{definition}

\begin{prop} \label{Majorization as convex hull of permutations}
For all $S\subseteq \mathfrak{a}^{d}_{+}$, we have that 
\begin{equation}
    \mathrm{Maj}(S)=\mathfrak{a}^{d}_{+}\cap \overline{\mathrm{conv}}(\mathrm{S}_{d}\cdot S),
\end{equation}
where $\overline{\mathrm{conv}}$ denotes the closed convex hull.
\end{prop}
\begin{proof}
    This follows immediately from Propositions \ref{MajorizationByConvexHull} and \ref{MconvexProposition}. 
\end{proof}

\begin{figure}[h]
\centering
\includegraphics[width=0.7\textwidth]{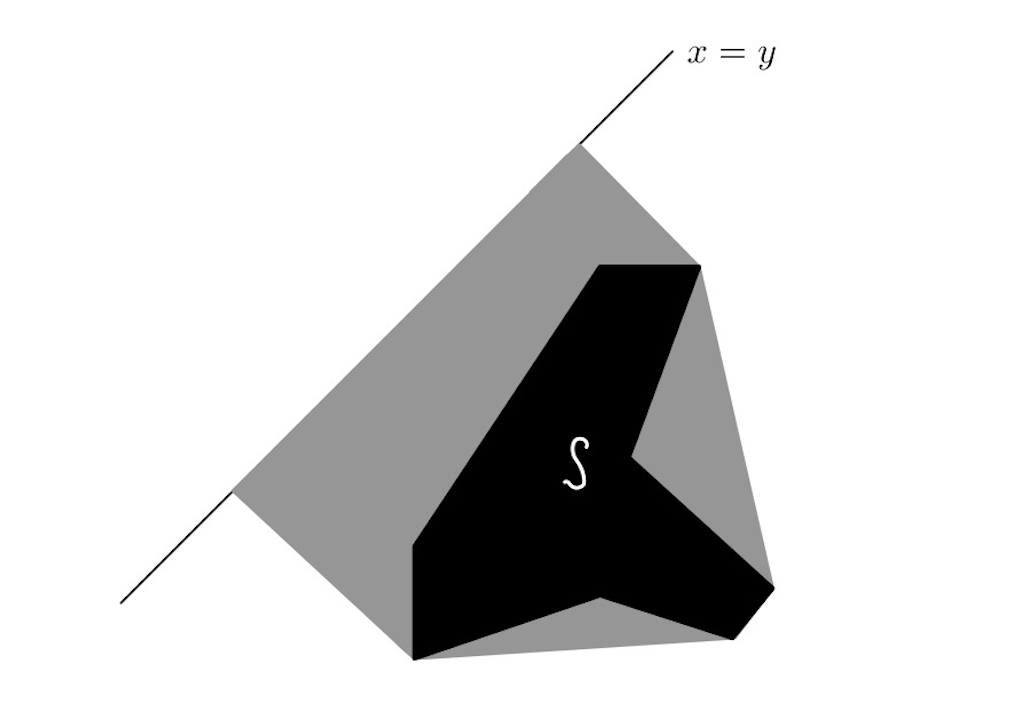}
\caption{Example of majorizarion hull on $\mathfrak{a}_{+}^{2}$.}
\label{Majorization hull two dimensional example}
\end{figure}

We illustrate this definition with two examples. In Figure \ref{Majorization hull two dimensional example}, we see an example of a majorization hull in two dimensions. The black silhouette is the original set $S\subseteq \mathfrak{a}_{+}^{2}$, while the gray set is the majorization hull $\mathrm{Maj}(S)$. Notice that in the two dimensional case, $\xi\in \mathrm{Maj}(S)$ if there exist $\eta\in\overline{\mathrm{conv}}(S)$ and $t\geq 0$ such that $\eta-\xi=t(1,-1)$.

\begin{figure}[h]
\centering
\includegraphics[width=0.7\textwidth]{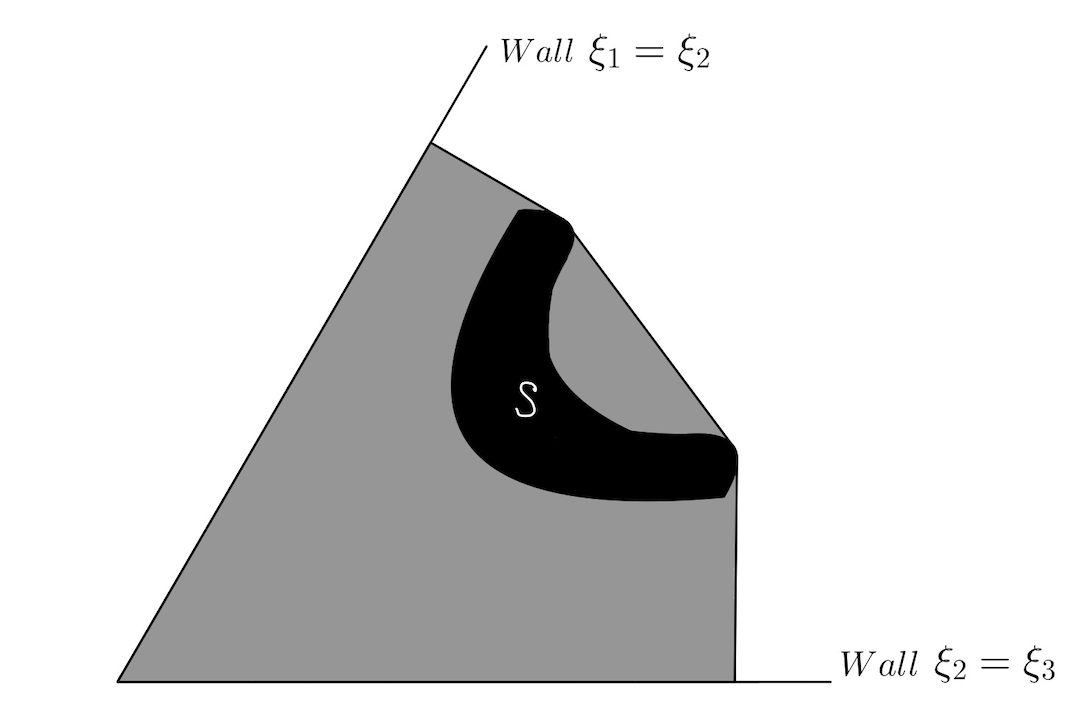}
\caption{Example of majorizarion hull on the fundamental Weyl chamber of $\mathfrak{sl}(3,\mathbb{R})$.}
\label{Majorization hull SL(3) example}
\end{figure}

Denote $\mathfrak{sl}(3,\mathbb{R})\coloneq \{(\xi_1,\xi_2,\xi_3)\in\mathbb{R}^{3}: \ \xi_1+\xi_2+\xi_3=0\}$. In Figure \ref{Majorization hull SL(3) example}, we see a set $S\subseteq \mathfrak{a}_{+}^{3}\cap\mathfrak{sl}(3,\mathbb{R})$, colored in black, and its majorization hull $\mathrm{Maj}(S)$ shaded in gray. The shape of $\mathrm{Maj}(S)$ can be determined using Proposition \ref{Majorization as convex hull of permutations}. 

Another characterization of the majorization hull will be given in Corollary \ref{MajorizationCharacterization}. 

The following result points to the existence of supersets of the uniform spectrum which may be simpler to study. 
\begin{prop}
Suppose that $T:X\to X$ is a continuous map and $\Phi:\mathbb{E}\to\mathbb{E}$ is a bundle map covering $T$. Then for all $n\in\mathbb{N}$ and any Riemannian metric $g$, the following inclusion holds
\begin{equation}
\mathcal{U}(\Phi)\subseteq \frac{1}{n}\mathrm{Maj}(\kappa^{g}(\Phi^n)).
\end{equation}
\end{prop}
(Compare with \cite[Lemma~2.2]{MisiurewiczZiemian}.)
\begin{proof}
Let $n$ be a natural number and let $v$ be a vector in the uniform spectrum. Then there exist $\{m_{j}\}\subseteq \mathbb{N}$ and $\{x_{m_{j}}\}\subseteq X$ such that $v=\lim_{j\to\infty}m_j^{-1}\kappa^{g}(\Phi^{m_j}_{x_{m_j}})$. Writing $m_j=nk_j+r_j$ for some integers $k_j$ and $0\leq r_j<n$ and using (\ref{TriangleInequalityCartanProjection}), we see that 

\begin{equation}
\dfrac{1}{m_j}\kappa^{g}(\Phi^{m_j}_{x_{m_j}})-\dfrac{1}{m_j}\kappa^{g}(\Phi^{r_j}_{T^{nk_j}x_{m_j}})\trianglelefteq\dfrac{nk_j}{m_j}\left(\dfrac{1}{n}\sum_{i=0}^{k_{j}-1}\dfrac{1}{k_j}\kappa^{g}(\Phi^{n}_{T^{in}x_{m_j}})\right).
\end{equation}
Note that the sequence 
\begin{equation}
\left\{\dfrac{nk_j}{m_j}\left(\dfrac{1}{n}\sum_{i=0}^{k_{j}-1}\dfrac{1}{k_j}\kappa^{g}(\Phi^{n}_{T^{in}x_{m_j}})\right)\right\}_{j\in\mathbb{N}}
\end{equation}
is bounded and thus has limit points. Hence, taking $j\to\infty$ and a subsequence if necessary
\begin{equation}
v\trianglelefteq \lim_{j\to\infty} \dfrac{1}{n}\sum_{i=0}^{k_{j}-1}\dfrac{1}{k_j}\kappa^{g}(\Phi^{n}_{T^{in}x_{m_j}}).
\end{equation}
The result follows from noticing that the right hand side belongs to $\frac{1}{n}\mathrm{Maj}(\kappa^{g}(\Phi^n))$.
\end{proof}
\subsection{Schur-convexity and separation} \label{Aura}

We will need to work with functionals $f:\mathfrak{a}_{+}^{d}\to\mathbb{R}$ which are \textit{monotonic} with the majorization relation. 

\begin{definition}
We will say that a function $f:\mathfrak{a}_{+}^{d}\to \mathbb{R}$ is \textit{Schur-convex} if 
\begin{equation}
  \xi\trianglelefteq\eta\implies f(\xi)\leq f(\eta).  
\end{equation}
\end{definition}

We are interested in Schur-convex functions which are also linear. For the rest of the paper, we will denote the regular inner product of $u,v\in\mathbb{R}^{d}$ as $u^{\top}v$. 

\begin{prop} \label{characterization of Schur-convex}
    A linear functional $f:\mathfrak{a}_{+}^{d}\to\mathbb{R}$ is Schur-convex if and only if there is ${b}\in \mathfrak{a}_{+}^{d}$ such that $f(\xi)={b}^{\top}\xi$.
\end{prop}
\begin{proof}
    Let $f:\mathfrak{a}_{+}^{d}\to\mathbb{R}$ be the restriction of a linear functional to $\mathfrak{a}_{+}^{d}$. Then there exist $b\in\mathbb{R}^{d}$ such that 
    \begin{equation}
        f(\xi)=b^{\top}\xi=\sum_{i=1}^{d}b_{i}\xi_{i}. 
    \end{equation}
    Define $a\in \mathbb{R}^{d}$ by 
    \begin{equation}
        a_{i}\coloneq\begin{cases}
            b_{i}-b_{i+1} &\text{ if }i<d,\\
            b_{d}&\text{ if }i=d.
        \end{cases}
    \end{equation}
    By the summation by parts formula, 
    \begin{equation}  \label{NewConstants}
        f(\xi)=\sum_{i=1}^{d}a_iv_i \ \text{ where } \ v_{i}\coloneq \sum_{j=1}^{i}\xi_{j}.
    \end{equation}

Let $P\subseteq \mathbb{R}^{d}$ be the set of vectors $v=(v_{1},\dots, v_{d})$ obtained as defined on (\ref{NewConstants}). This set is isomorphic to $\mathfrak{a}_{+}^{d}.$

Take a vector $\eta\in\mathfrak{a}_{+}^{d}$ and let $w_{i}\coloneq \sum_{j=1}^{i}\eta_{j}$. By definition, $\xi\trianglelefteq\eta$ if and only if $v_{i}\leq w_{i}$ for each $i=1,\dots, d-1$ and $v_{d}=w_{d}$.

Therefore, $f$ is Schur-convex if and only if the function $v\in P\mapsto \sum_{i=1}^{d}a_{i}v_{i}$ is increasing with respect to the variables $v_{1},\dots, v_{d-1}$. This condition is equivalent to $a_{1},\dots, a_{d-1}\geq 0$, that is, $b\in \mathfrak{a}_{+}^{d}.$
\end{proof}

\begin{remark} 
    The representation (\ref{NewConstants}) is often convenient. The coefficients $a_{i}$ are given directly by the formula $a_{i}=-f(\hat{e}_{i-1})+2f(\hat{e}_{i})-f(\hat{e}_{i+1})$, where $\hat{e}_{i}\coloneq (\underbrace{1,\dots,1}_{i}, \underbrace{0,\dots, 0}_{d-i})\in \mathfrak{a}_{+}^{+}$ and $\hat{e}_{d+1}\coloneq \hat{e}_{d}$. 
    
\end{remark}

The following proposition and corollary are analogues of standard results in convex analysis \cite[p.~49]{BoydV} to the setting of LCM-closed sets. 
 
\begin{prop}[Separation theorem] \label{SeparationTheorem}
    Let $C\subseteq \mathfrak{a}_{+}^{d}$ be a nonempty LCM-closed set and let $\xi\in\mathfrak{a}_{+}^{d}\backslash C$ . Then there exist a linear Schur-convex function $f$ such that $f(\xi)>\sup_{C}f.$
\end{prop}
\begin{proof}
By Proposition~\ref{MconvexProposition}, the set $\hat{C} \coloneq \mathrm{S}_d \cdot C$ is closed and convex. 
Let $\eta$ be the unique point in $\hat{C}$ that minimizes the Euclidean distance $\|\xi - \eta\|$.
We claim that $\eta \in \mathfrak{a}_{+}^{d}$ and, in particular, $\eta \in C$.
Fix $k \in \{1,\dots,d-1\}$; we want to check that $\eta_k \ge \eta_{k+1}$.
Let $\hat \eta$ be obtained from the vector $\eta$ by exchanging the $k$\textsuperscript{th} and $k+1$\textsuperscript{th} components; then $\hat \eta \in \hat{C}$.
If $\eta_k \neq \eta_{k+1}$, then $\hat{\eta} \neq \eta$ and so
\begin{align}
\begin{split}
0
&< \|\xi - \hat{\eta}\|^2 - \|\xi - \eta\|^2  \\ 
&= (\xi_k-\eta_{k+1})^2 + (\xi_{k+1}-\eta_k)^2 - (\xi_k-\eta_k)^2 - (\xi_{k+1}-\eta_{k+1})^2 \\
&= 2(\xi_k-\xi_{k+1})(\eta_k - \eta_{k+1}) \, .
\end{split}
\end{align}	
In particular, $\eta_k > \eta_{k+1}$, completing the proof of the claim. 

Define the linear functional $f(\zeta) \coloneq (\xi-\eta)^\top \zeta$.
Note that $\sup_C f$ is attained at the unique point $\eta$, and that $f(\eta)<f(\xi)$.
To complete the proof, we will show that $f$ is Schur-convex, which by Proposition~\ref{characterization of Schur-convex} means that $\xi - \eta \in \mathfrak{a}_{+}^{d}$.
Fix $k \in \{1,\dots,d-1\}$; we want to check that $\xi_k-\eta_k \ge \xi_{k+1}-\eta_{k+1}$.
This inequality is satisfied if $\eta_k = \eta_{k+1}$, so assume that $\eta_k > \eta_{k+1}$.
Define $\hat \eta$ exactly as above. 
Then 
\begin{align}
\begin{split}
0
&< f(\eta) - f(\hat \eta) \\
&= (\xi-\eta)^\top (\eta - \hat \eta) \\ 
&= (\xi_k-\eta_k)(\eta_k-\eta_{k+1}) + (\xi_{k+1}-\eta_{k+1})(\eta_{k+1}-\eta_k)  \\
&= (\xi_k-\eta_k-\xi_{k+1}+\eta_{k+1}) (\eta_k-\eta_{k+1}) \, , 
\end{split}
\end{align}	
and the desired inequality follows.
\end{proof}

From the separation theorem, we get the following characterization for the majorization hull.  

\begin{corollary}[Half-space representation] \label{MajorizationCharacterization}
    Given a set $S\subseteq \mathfrak{a}_{+}^{d}$ we have
    \begin{equation}
    \mathrm{Maj}(S)=\bigcap_{\substack{f\text{ Schur-convex} \\ \text{linear} }}\left\{\xi\in\mathfrak{a}_{+}^{d}: f(\xi)\leq \sup_{S}f\right\}.    
    \end{equation}
\end{corollary}
\begin{proof}
    For each Schur-convex linear functional $f$, the set $\{\xi \in \mathfrak{a}_{+}^{d} : f(\xi) \le \sup_S f \}$ is LCM-closed and contains $S$. Let $C$ be intersection of this family of sets; then $C$ is LCM-closed and contains $S$. By definition, $\mathrm{Maj}(S)$ is the smallest set with such properties; therefore, $\mathrm{Maj}(S) \subseteq C$. On the other hand, suppose that $\xi \not\in \mathrm{Maj}(S)$. By Proposition~\ref{SeparationTheorem}, there exists a linear Schur-convex functional $f$ such that $f(\xi) > \sup_{\mathrm{Maj}(S)} f \ge \sup_S f$.
That is, $\xi \not\in C$. Therefore, $\mathrm{Maj}(S)=C$.
\end{proof}

\subsection{Bounding the uniform spectrum in terms of the Lyapunov spectrum}

Assume that $T:X\to X$ is a continuous map, $\Phi:\mathbb{E}\to\mathbb{E}$ is a bundle map and $f:\mathfrak{a}_{+}^{d}\to\mathbb{R}^{d}$ is a linear Schur-convex function. 

Take $n,m\in \mathbb{N}$. We see that for every $x\in X$, $\Phi_{x}^{m+n}=\Phi_{T^{n}x}^{m}\circ \Phi^{n}_x$. By (\ref{TriangleInequalityCartanProjection})
\begin{equation}
    \kappa^{g}(\Phi_{x}^{m+n})\trianglelefteq \kappa^{g}(\Phi_{T^{n}x}^{m})+ \kappa^{g}(\Phi^{n}_x).
\end{equation}
Thus, by Schur-convexity
\begin{equation}
    f(\kappa^{g}(\Phi_{x}^{m+n}))\leq f(\kappa^{g}(\Phi_{T^{n}x}^{m}))+ f(\kappa^{g}(\Phi^{n}_x)).
\end{equation}
Hence, we conclude that the sequence of continuous real-valued functions 
\begin{equation}\label{Iteratedfunctionals}
    F_{n}(x)\coloneq f(\kappa^{g}(\Phi^{n}_x)) 
\end{equation}
is subbadditive with respect to the dynamics coming from $T$.

Define 
\begin{equation} \label{definition of top}
    \lambda^{\top}_{f}\coloneq \lim_{n\to\infty}\dfrac{1}{n}\sup_{x\in X}f(\kappa^{g}(\Phi^{n}_x)).
\end{equation}
The existence of the limit is ensured by Fekete's lemma. 
\begin{prop} \label{proposition about top}
Let $T:X\to X$ be a continuous map, $\Phi:\mathbb{E}\to\mathbb{E}$ be a bundle map, and $f:\mathfrak{a}_{+}^{d}\to\mathbb{R}$ be a linear Schur-convex function. Then
\begin{equation} \label{top}
    \lambda_{f}^{\top}=\max_{\mathcal{L}(\Phi)}f= \max_{\mathcal{U}(\Phi)}f.
\end{equation}
\end{prop}
\begin{proof}
    First, recall that from the Furstenberg-Kesten theorem there exists $\varphi(\mu, \Phi)\in \mathfrak{a}_{+}^{d}$ such that 
    \begin{equation} \label{Lyapunov vector}
    \varphi(\mu, \Phi)=\lim_{n\to\infty}\dfrac{1}{n}\int\kappa^{g}(\Phi^{n}_{x})d\mu.
    \end{equation}
    
    Additionally, given any subbadditive sequence $\{F_n\}$ of upper semicontinuous functions, \cite[Theorem~A.3]{MorrisErgodic} states that 
    \begin{equation}
        \lim_{n\to\infty}\dfrac{1}{n}\sup_{X}F_{n}=\max_{\mu\in \mathcal{M}_{\mathrm{erg}}(T)} \lim_{n\to\infty}\dfrac{1}{n}\int F_{n}d\mu.
    \end{equation}
    We specialize this result to the sequence defined in (\ref{Iteratedfunctionals}). Note that for any $\mu\in \mathcal{M}_{\mathrm{erg}}(T)$, by (\ref{Lyapunov vector})
    \begin{equation}
        \lim_{n\to\infty}\dfrac{1}{n}\int F_{n}d\mu=f(\varphi(\mu,\Phi)). 
    \end{equation}
    Therefore, $\lambda_{f}^{\top}=\max_{\mathcal{L}(\Phi)}f$, proving the first equality in (\ref{top}).

   Since $\mathcal{U}(\Phi)$ is a compact set which contains $\mathcal{L}(\Phi)$, $\max_{\mathcal{U}(\Phi)}f$ is well defined and at least $\lambda_{f}^{\top}=\max_{\mathcal{L}(\Phi)}f$. 

   Pick a point $v\in\mathcal{U}(\Phi)$ where $f$ attains the maximum. Then there are sequences $\{x_{j}\}\subseteq X$ and $\{n_{j}\}\subseteq \mathbb{N}$ where $n_{j}\to\infty$ such that $v=\lim_{j\to\infty}\frac{1}{n_j}\kappa^{g}(\Phi^{n_j}_{x_j})$. Then
   \begin{equation}
       \lambda^{\top}_{f}\leq f(v)=\lim_{j\to\infty}\frac{1}{n_j}f(\kappa^{g}(\Phi^{n_j}_{x_j}))\leq \lim_{j\to\infty}\frac{1}{n_j}\sup_{x\in X}f(\kappa^{g}(\Phi^{n_j}_{x}))=\lambda_{f}^{\top},
   \end{equation}
   completing the proof. 
\end{proof}
\begin{corollary} \label{First bound of uniform spectrum}
    The uniform spectrum is contained in the majorization hull of the Lyapunov spectrum. That is,
    \begin{equation}
    \mathcal{U}(\Phi)\subseteq \mathrm{Maj}(\mathcal{L}(\Phi)).
    \end{equation}
    In particular, since $\mathcal{L}(\Phi)\subseteq \mathcal{U}(\Phi)$, we have 
    \begin{equation}
        \mathrm{Maj}(\mathcal{U}(\Phi))=\mathrm{Maj}(\mathcal{L}(\Phi)).
    \end{equation}
\end{corollary}
(Compare with \cite[Theorem~2.4]{MisiurewiczZiemian}.)
\begin{proof}
    From Proposition \ref{proposition about top}
    \begin{align}
    \begin{split}
        \mathcal{U}(\Phi)&\subseteq\bigcap_{\substack{f\text{ Schur-convex} \\ \text{linear} }}\left\{\xi\in\mathfrak{a}_{+}^{d}: f(\xi)\leq \max_{\mathcal{U}(\Phi)}f\right\}\\
        &=\bigcap_{\substack{f\text{ Schur-convex} \\ \text{linear} }}\left\{\xi\in\mathfrak{a}_{+}^{d}: f(\xi)\leq \max_{\mathcal{L}(\Phi)}f\right\}\\
        &=\mathrm{Maj}(\mathcal{L}(\Phi)),
    \end{split}
    \end{align}
    where the last equality comes from the characterization given in Proposition \ref{MajorizationCharacterization}.
\end{proof}
We end this subsection stating the following result from Bochi \cite[Proposition~4.1]{BochiErgodicOptimization} which gives another relation between the uniform spectrum and the Cartan projection of iterates of the bundle map. 

\begin{prop} \label{Change of norm}
    Assume that $T:X\to X$ is a continuous map and $\Phi:\mathbb{E}\to\mathbb{E}$ is a bundle map covering $T$. Given any open set $\mathcal{O}\subseteq \mathfrak{a}_{+}^{d}$ such that $\mathrm{Maj}(\mathcal{U}(\Phi))\subseteq \mathcal{O}$, the following two conclusions hold:
    \begin{enumerate}
        \item There exist $N\in\mathbb{N}$ such that $n^{-1}\kappa^{g}(\Phi^{n})\subseteq \mathcal{O}$ for all $n\geq N$.
        \item There exist a $C^{0}$ Riemannian metric $\tilde{g}$ such that $\kappa^{\tilde{g}}(\Phi)\subseteq \mathcal{O}$.
    \end{enumerate}
\end{prop}
(Compare part (1) with \cite[Lemma~2.8]{MisiurewiczZiemian} and part (2) with \cite[Proposition~4.2]{Schub}).

\subsection{Aura}
We will need the following definition for the next sections. 

\begin{definition} \label{definition of aura}
    Let $S \subset \mathfrak{a}_{+}^{d}$ be a bounded set and $\gamma\geq 1$. The $\gamma$-aura of $S$ is the set
    \begin{equation}
    \hspace{-0.3cm}\mathrm{Aura}_{\gamma}(S)\coloneq\bigcap_{\substack{f\text{ Schur-convex} \\ \text{linear} }}\left\{\xi\in\mathfrak{a}_{+}^{d}: (1-\gamma)\sup_{S}f+\gamma\inf_{S}f\leq f(\xi)\leq \sup_{S}f\right\}.
    \end{equation}
\end{definition}

Observe the following properties:
\begin{itemize}
    \item For any $1\leq \gamma\leq \gamma'$ we obtain that $\text{Aura}_{1}(S)\subseteq \text{Aura}_{\gamma}(S)\subseteq \text{Aura}_{\gamma'}(S).$
    \item From Corollary \ref{MajorizationCharacterization} it is clear that $\text{Aura}_{\gamma}(S)\subseteq \mathrm{Maj}(S)$ for all $\gamma\geq 1$.
    \item The $\gamma-$aura of $S$ is always a compact and convex set, which readily implies that $\overline{\text{conv}}(S)\subseteq \text{Aura}_{1}(S).$
\end{itemize}

Let us illustrate the definition with a few examples. 

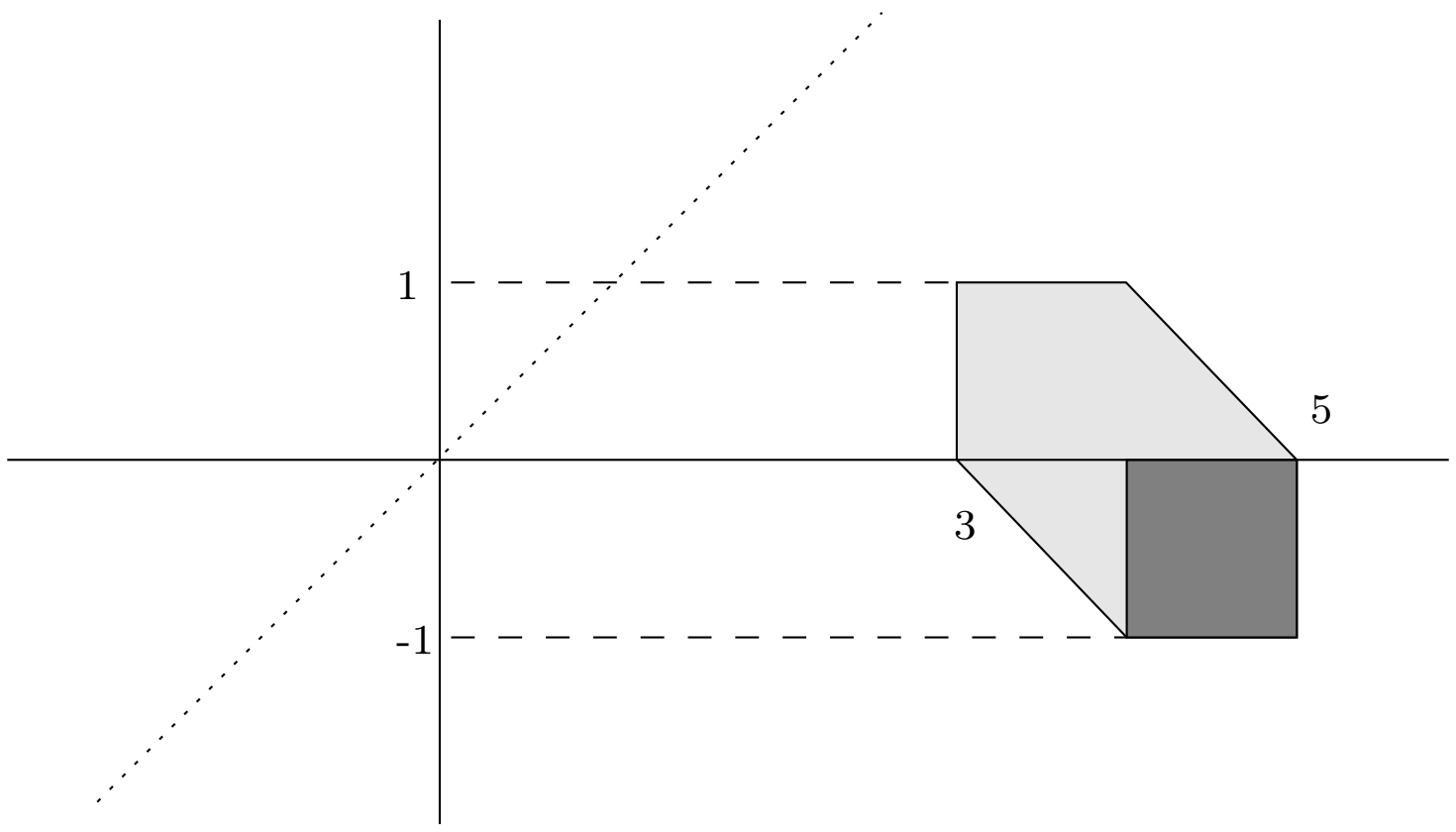
\begin{figure}[h] 

\centering

\begin{tikzpicture}[x=0.75pt,y=0.75pt,yscale=-1,xscale=1]


\draw[thin, fill=gray!20] (529.35,176.21) -- (529.35,221.19) -- (486.27,221.19) -- (443.2,176.21) -- (443.2,131.27) -- (486.09,131.27) -- cycle;

\draw[thin, fill=gray] (486.27,176.21) -- (529.35,176.21) -- (529.35,221.19) -- (486.27,221.19) -- cycle ;

\draw  [dash pattern={on 4.5pt off 4.5pt}]  (315.09,131.27) -- (443.2,131.27) ;
\draw  [dash pattern={on 4.5pt off 4.5pt}]  (315.09,221.19) -- (486.27,221.19) ;

\draw (300,127.27) node [anchor=north west][inner sep=0.75pt]  [font=\footnotesize] [align=left] {1};
\draw (300,217.19) node [anchor=north west][inner sep=0.75pt]  [font=\footnotesize] [align=left] {\mbox{-}1};
\draw (441.28,188.06) node [anchor=north west][inner sep=0.75pt]  [font=\footnotesize] [align=left] {3};
\draw (531.46,158.82) node [anchor=north west][inner sep=0.75pt]  [font=\footnotesize] [align=left] {5};

\draw [thin]   (312.22,64.71) -- (312.22,268.5) ;

\draw [thin]   (202.64,176.21) -- (567.86,176.21) ;

\draw  [dash pattern={on 0.84pt off 2.51pt}]  (225.46,262.88) -- (424.28,62.91) ;

\end{tikzpicture}
\caption{The figure shows $S$ in dark gray and Aura$_{2}(S)$ in light gray. Note that to calculate $\mathrm{Aura}_{2}(S)$ it is only necessary to evaluate the Schur-convex linear functionals at the boundary of $S$.}
\label{SwellingExample}
\end{figure}

For a two-dimensional example, consider the square $S=[4,5]\times [-1,0]\subseteq \mathfrak{a}_{+}^{d}.$ and $\gamma=2$. Then the set $\text{Aura}_{2}(S)$ is the hexagon shown in Figure \ref{SwellingExample}.

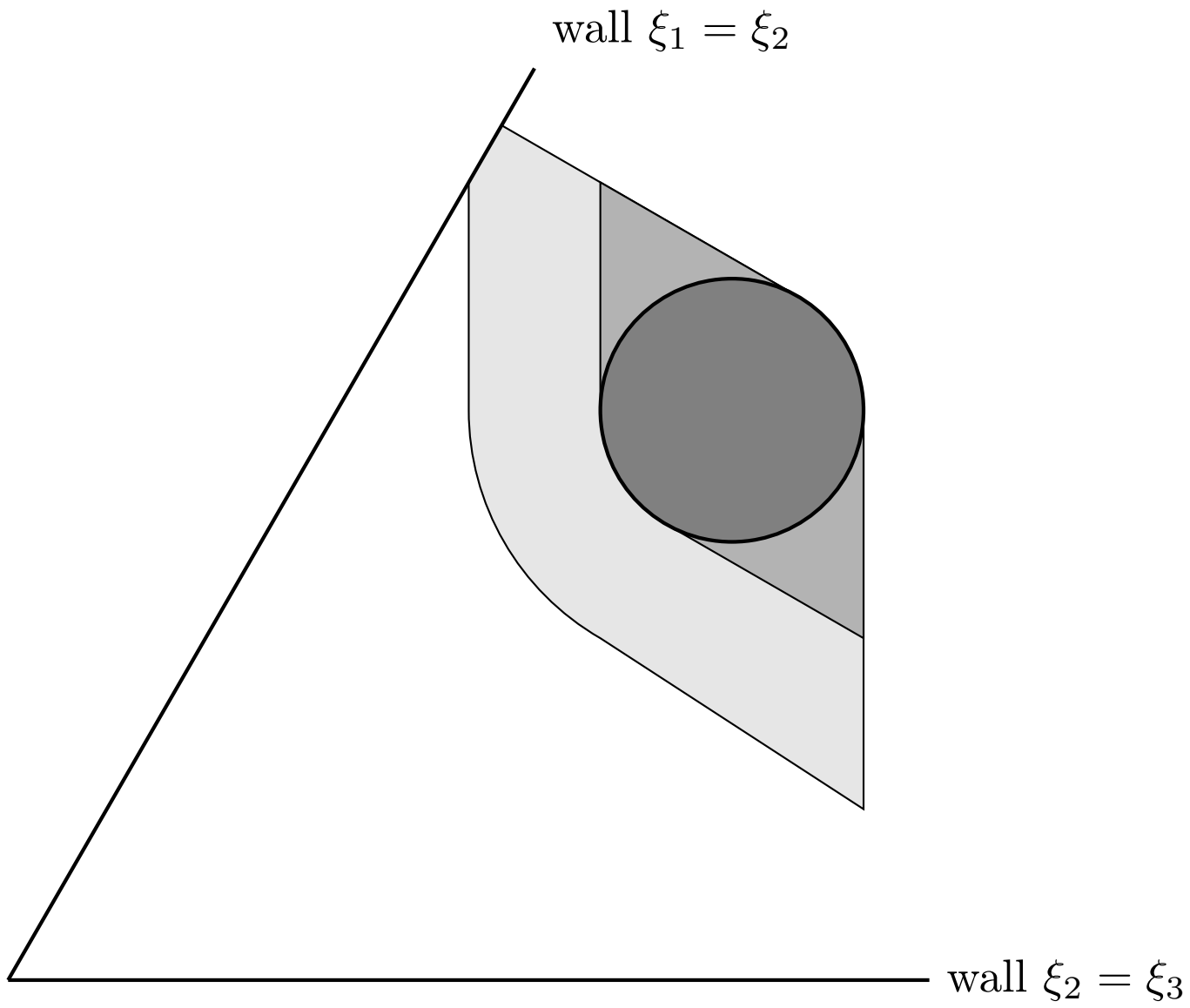
\begin{figure} 
	\begin{tikzpicture}
		\def\sqrtthree{1.732}

		
		\draw[fill=gray!20] (1,0) arc (0:60:1) -- (-1.75,1.25*\sqrtthree) -- (-2,\sqrtthree) -- (-2,0) arc (180:240:2) -- (1,-1.75*\sqrtthree) -- cycle;

		\draw[fill=gray!60] (1,0) arc (0:60:1) -- (-1,\sqrtthree) -- (-1,0) arc (180:240:1) -- (1,-\sqrtthree) -- cycle;

		\draw[thick, fill=gray] (0,0) circle (1);

		\draw[thick] (-5.5,-2.5*\sqrtthree) --++ (7,0) node[right]{wall $\xi_2=\xi_3$};
		\draw[thick] (-5.5,-2.5*\sqrtthree) --++ (4,4*\sqrtthree) node[above right]{wall  $\xi_1=\xi_2$};

	
	\end{tikzpicture}
	\caption{Letting $d=3$ and $S$ be a round disk (pictured in dark gray) in the plane, the figure shows $\mathrm{Aura}_1(S)$ (light gray) and $\mathrm{Aura}_{1.5}(S)$ (lighter gray).}
    \label{Eye}
\end{figure}

An example where $S \subseteq \mathfrak{a}_{+}^{3}\cap \mathfrak{sl}(3,\mathbb{R})$ is shown in Figure \ref{Eye}. 

In the case that $S$ lives in a one-dimensional or two-dimensional space, it is possible to easily characterize $\mathrm{Aura}_1(S)$, as the following proposition shows.  

\begin{prop}\label{Aura in dimension 2}
    Suppose that $d=1,2$ and let $S\subseteq \mathfrak{a}_{+}^{d}$ be a bounded set. Then $\mathrm{Aura}_{1}(S)=\overline{\mathrm{conv}}(S)$. 
\end{prop}
\begin{proof}
 In case that $d=1$, all the linear functions are Schur-convex, and by the characterization of the convex closure as the intersection of half-spaces which contain $S$, we obtain that $\mathrm{Aura}_{1}(S)\subseteq\overline{\mathrm{conv}}(S).$ 
 
 When $d=2$, we see that if a linear function $f$ is not Schur-convex, then $-f$ is Schur-convex. Thus
    \begin{align}
    \begin{split}
        \overline{\mathrm{conv}}(S)&=\bigcap_{f\text{ linear}}\left\{\xi\in\mathfrak{a}_{+}^{d}: f(\xi)\leq \sup_{S}f\right\}\\
        &=\left(\bigcap_{\substack{f\text{ Schur-convex} \\ \text{linear} }}\left\{\xi\in\mathfrak{a}_{+}^{d}: f(\xi)\leq \sup_{S}f\right\}\right)\cap\\
        &\hspace{2cm} \left(\bigcap_{\substack{-f\text{ Schur-convex} \\ \text{linear} }}\left\{\xi\in\mathfrak{a}_{+}^{d}: f(\xi)\leq \sup_{S}f\right\}\right)
        \end{split}\\
        &=\bigcap_{\substack{f\text{ Schur-convex} \\ \text{linear} }}\left\{\xi\in\mathfrak{a}_{+}^{d}: \inf_{S}f\leq f(\xi)\leq \sup_{S}f\right\}=\mathrm{Aura}_{1}(S). \qedhere
    \end{align}
\end{proof}
\begin{remark}
    Proposition \ref{Aura in dimension 2}, fails in higher dimensions. This is clear from Figure \ref{Eye}. For another example, take $\xi=(3,2,1)$, $\omega=(3,3,0)$ and $\eta=(4,2,0)$. Note that $\xi\trianglelefteq\omega\trianglelefteq\eta$ and thus $\omega\in \mathrm{Aura}_{1}(\{\xi,\eta\})$. However, note also that $\omega\not\in \mathrm{conv}(\{\xi,\eta\}).$
\end{remark}

\section{Improved bounds for the uniform spectrum} \label{Section 4}

Recall that $X$ is a compact metric space. An algebra of continuous real-valued functions that contains the constant 
functions is called \textit{normal} if it satisfies Urysohn lemma and consequently admits partitions of unity \cite[page~221]{Katznelson}. 

\textit{Standing hypothesis}:  $\theta>0$ is a number so that the algebra of $\theta$-H\"older functions is normal. 
This hypothesis is automatically satisfied if either $0<\theta\leq 1$ or $X$ has zero topological dimension.

\subsection{Hyperbolic homeomorphisms of $\theta$-H\"older vector bundles}

Consider a finite family of charts $\{U_i,\psi_i\}$ of $X$ and for all $x\in X$ recall the isomorphisms $h_{i}:\mathbb{R}^{d}\to \mathbb{E}_{x}$ given in (\ref{definition of h}). We define $\mathbb{E}$ as a $\theta$-\textit{Hölder vector bundle} if the \textit{coordinate transformations}
\begin{equation}
g_{j\xleftarrow{}i}(x)=h_{j}(x)^{-1}\circ h_{i}(x):U_j\cap U_i\to \text{GL}(d,\mathbb{R}),
\end{equation}
are $\theta$-H\"older maps. 

Furthermore, we call a bundle map $\Phi$ as \textit{$\theta$-Hölder} if $T$ is Lipschitz and the maps 
\begin{equation}
x\mapsto h_{j}(Tx)^{-1}\circ \Phi_x \circ h_{i}(x)
\end{equation}
are $\theta$-Hölder.

Given $\epsilon>0$ and $x\in X$, we define the local stable set and the local unstable set respectively as
\begin{align}
\begin{split}
    W^{s}_{\epsilon}(x)&\coloneq\{y\in X:d(T^{n}y,T^{n}x)\leq\epsilon \text{ for all }n\geq 0\},\\
    W^{u}_{\epsilon}(x)&\coloneq\{y\in X:d(T^{-n}y,T^{-n}x)\leq\epsilon \text{ for all }n\geq 0\}.
\end{split}
\end{align}
\begin{definition}
    We will call $T$ a \textit{hyperbolic homeomorphism } if the following hold:
    \begin{enumerate}
        \item $T$ is bi-Lipschitz. 
        \item There exist $0<\epsilon_0<1$ and $\tau>0$ such that for all $x\in X$
        \begin{align}
        \begin{split}
            w,z\in W^{s}_{\epsilon_0}(x)&\implies d(Tw, Tz)\leq e^{-\tau}d(w,z),\\
            w,z\in W^{u}_{\epsilon_0}(x)&\implies d(T^{-1}w, T^{-1}z)\leq e^{-\tau}d(w,z).
        \end{split}
        \end{align}
        \item There exist a constant $\epsilon_{1}\in (0,\epsilon_{0})$ such that for any $x,y\in X$ with $d(x,y)\leq 2\epsilon_1$, there is a unique element in the set $W^{s}_{\epsilon_{0}}(x)\cap W^{u}_{\epsilon_{0}}(y)$ which is denoted $[y,x]$. Moreover, the bracket function $[\cdot,\cdot]$ is continuous on both entries. 
        \item There is a constant $C>0$ such that for $x,y\in X$ where $d(x,y)\leq 2\epsilon_1$ then 
        \begin{equation}
        \max\{d([x,y],x),d([x,y],y)\}\leq Cd(x,y).
        \end{equation}
    \end{enumerate}
\end{definition}

The key property of hyperbolic homeomorphisms that we will need is the following exponential closing lemma. 

\begin{lemma} \label{approximation}
    Suppose $(X,d)$ is a compact metric space and $T:X\to X$ is a topologically mixing hyperbolic homeomorphism. Then there are $C>0$ and $n_{0}\in \mathbb{N}$ such that for $n\geq n_{0}$ and all $x\in X$ we can find a periodic point $p=p(x)$ of period $n$ such that for $0\leq i\leq n$
    \begin{equation}
      d(T^ix, T^ip)\leq Ce^{-\tau\min\{i,n-i\}}.
      \end{equation}
\end{lemma}

We were not able to find an exact reference for this lemma, so we provide a proof in Appendix \ref{Appendix B}. 

\subsection{Statement of the theorem}

Note that the linear functional $g(\xi_{1},\dots, \xi_{d})=\xi_{1}-\xi_{d}$ is Schur-convex. In terms of Proposition \ref{characterization of Schur-convex}, it is represented by the vector $b=(1,0,\dots, 0,-1)$. Denote $(\xi_{1}-\xi_{d})^{\top}\coloneq \max_{\mathcal{U}(\Phi)}g$. By Proposition \ref{proposition about top}, this quantity equals $\lambda_{g}^{\top}$. 
\begin{theorem} \label{MainTheorem}
    Assume that $T:X\to X$ is a topologically transitive hyperbolic homeomorphism with constant of expansion $\tau>0$, $\Phi:\mathbb{E}\to\mathbb{E}$ is a $\theta$-Hölder bundle map covering $T$ with $\theta>0 $ and $\gamma_{0}=\max\{1,(\xi_1-\xi_d)^{\top}/\tau\theta\}$. Then,
    \begin{equation}
    \mathcal{U}(\Phi)\subseteq \mathrm{Aura}_{\gamma_{0}}(\mathcal{L}(\Phi)).
    \end{equation}
\end{theorem}

The Theorems \cite[Theorem~3.10]{DewittGogolev} from DeWitt and Gogoglev and \cite[Theorem~1]{Renato} from Velozo Ruiz can be obtained as corollaries of the previous theorem as we will explain in Section \ref{Corollaries}.  

Let us reformulate Theorem \ref{MainTheorem} in terms of the supremum and infimum values of linear Schur-convex functions evaluated in $\mathcal{L}(\Phi)$ and $\mathcal{U}(\Phi)$. Recall from Proposition \ref{proposition about top} that under the conditions of Theorem \ref{MainTheorem}
\begin{equation}
     \lambda_{f}^{\top}=\max_{\mathcal{L}(\Phi)}f= \max_{\mathcal{U}(\Phi)}f,
\end{equation}
for $f$ a linear Schur-convex function. We will also write
\begin{equation}
    \lambda_{f}^{\bot}=\inf_{\mathcal{L}(\Phi)}f \ \text{ and } \ \alpha_{f}^{\bot}=\min_{\mathcal{U}(\Phi)}f. 
\end{equation}
Unlike the case of $\lambda^{\top}_{f}$, it is possible that $\lambda_{f}^{\bot}$ is not equal to $\alpha_{f}^{\bot}$. Hence, in general we only know that $\alpha_{f}^{\bot}\leq \lambda_{f}^{\bot}$. 

Consider $\xi\in \mathcal{U}(\Phi)$. Then, for any $f$ linear Schur-convex function $f$, we have $\alpha_{f}^{\bot}\leq f(\xi) \leq \lambda^{\top}_{f}$.
Hence, $\mathcal{U}(\Phi)\subseteq \mathrm{Aura}_{\gamma}(\mathcal{L}(\Phi))$ if and only if
\begin{equation}
(1-\gamma)\lambda^{\top}_{f}+\gamma\lambda_{f}^{\bot}\leq \alpha_{f}^{\bot}\leq f(\xi).
\end{equation}

Furthermore, by the spectral decomposition theorem for hyperbolic homeomorphisms \cite[page~101]{AokiHiraide}, there is $k\in\mathbb{N}$ and a partition $X=X_{0}\sqcup\dots X_{k-1}$ such that $T(X_{i})=X_{i+1}$ modulo $k$ for all $i\in \{0,\dots, k-1\}$ and $T^{k}_{|X_{0}}:X_{0}\to X_{0}$ is topologically mixing. Note that $\nu\in \mathcal{M}_{\mathrm{erg}}(T^{k}_{|X_{0}})$ if and only if $T^{i}_{*}\nu\in \mathcal{M}_{\mathrm{erg}}(T^{k}_{|X_{i}})$ for each $i\in\{1,\dots,k-1\}$. Hence, for any $0\leq i,j\leq k-1$ we obtain that $\mathcal{L}(\Phi^{k}_{|\pi^{-1}(X_{i})})=\mathcal{L}(\Phi^{k}_{|\pi^{-1}(X_{j})})$ and therefore, from Proposition \ref{change of uniform spectrum due to iterations}, we see that $\mathcal{L}(\Phi^{k}_{|\pi^{-1}(X_{0})})=\mathcal{L}(\Phi^{k})=k\mathcal{L}(\Phi)$. On the other hand, from (\ref{TriangleInequalityCartanProjection}) and Proposition \ref{change of uniform spectrum due to iterations}, we can also conclude that $\mathcal{U}(\Phi^{k}_{|\pi^{-1}(X_{0})})=\mathcal{U}(\Phi^{k})=k\mathcal{U}(\Phi)$. Then, instead of the original system, we can study the bundle map $\Phi_{|\pi^{-1}(X_{0})}:\pi^{-1}(X_{0})\to \pi^{-1}(X_0)$ over the topologically mixing hyperbolic homeomorphism $T^{k}_{|X_{0}}:X_{0}\to X_{0}$ whose constant of expansion is $k\tau>0$. That is, to prove Theorem \ref{MainTheorem}, it is enough to prove the following result.

\begin{theorem} \label{Reformulation of main theorem}
    Assume that $T:X\to X$ is a topologically mixing hyperbolic homeomorphism and $\Phi:\mathbb{E}\to\mathbb{E}$ is a $\theta$-Hölder bundle map covering $T$. Let $f:\mathfrak{a}_{+}^{d}\to\mathbb{R}$ be a Schur-convex linear functional. Then, for all $\gamma\in \mathbb{R}$ such that
    \begin{equation}\gamma>\max\left\{1, \frac{(\xi_1-\xi_d)^{\top}}{\tau\theta}\right\},
    \end{equation}
    we obtain that
    \begin{equation}
    \lambda^{\top}_{f}-\alpha_{f}^{\bot}\leq \gamma(\lambda^{\top}_{f}-\lambda_{f}^{\bot}).
    \end{equation}
\end{theorem}

The proof of Theorem \ref{Reformulation of main theorem} will take the remainder of this section. 
\subsection{Growth of the bundle map restricted to vector subspaces}

The use of linear Schur-convex functional will allow us to understand the growth of the bundle map $\Phi$ along the directions of the singular values.

Let $f:\mathfrak{a}_{+}^{d}\to\mathbb{R}^{d}$ be a linear Schur-convex function. Recall the subbadditive sequence defined in (\ref{Iteratedfunctionals}):

\begin{equation}
    F_{n}(x)= f(\kappa^{g}(\Phi^{n}_x))=b^{\top}\kappa^{g}(\Phi^{n}_x). 
\end{equation}
Recall also the constants $a_{1},\dots, a_{d}\in\mathbb{R}$, obtained in (\ref{NewConstants}). Using these constants, we obtain the sequence
\begin{align} \label{F in terms of wedge}
        F_{n}(x)=\sum_{i=1}^{d}a_i\norm{\mathsf{\Lambda}^{i}\Phi^{n}_x},
\end{align}
where the norm of $\mathsf{{\Lambda}^{i}}\mathbb{E}_{x}$ is given by the inner product coming from
\begin{equation}
    \langle v_{1}\wedge\dots \wedge v_{i}, w_{1}\wedge\dots\wedge w_{i}\rangle=\det(g_x(v_i,w_j)),
\end{equation}
where $g$ is the Riemannian norm on $\mathbb{E}_{x}$. For more information, see \cite[Chapter~III]{Yokonuma}. 

This sequence will allow us to evaluate how the linear Schur-convex function $f$ changes with successive iterations of the bundle map. Our argument, however, also requires understanding the behavior $f$ after successive iterations of $\Phi$ along particular vector subspaces. Hence, we give the following definitions. 

\begin{definition}
    Let $x\in X$. A \textit{complete flag} of $\mathbb{E}_{x}$ is a tuple $\mathcal{V}=(V^{i})_{i\in \{0,\dots, d\}}$ where $V^{0},\dots, V^{d}\subseteq \mathbb{E}_{x}$ are vector subspaces such that $\dim V^{i}=i$ for $i\in\{0,\dots, d\}$ and $V^{i}\subseteq V^{i+1}$.
    We will write $\mathfrak{F}_{x}$ for the spaces of all complete flags on $\mathbb{E}_{x}$. 
\end{definition}

\begin{definition}
    Let $f:\mathfrak{a}_{+}^{d}\to\mathbb{R}$ be a linear Schur-convex functional and let $a_1,\dots, a_d\in \mathbb{R}$ be the constants defined on  (\ref{NewConstants}) and take $x,y\in X$. We define the \textit{induced function of $f$ on complete flags} as 
    \begin{align}
    \begin{split}
        \tilde{f}:\text{GL}(\mathbb{E}_x,\mathbb{E}_{y})\times \mathfrak{F}_{x}&\to \mathbb{R}\\
        (A,\mathcal{V})&\mapsto\sum_{i=1}^{d}a_i\log|\det(A|V^{i})|,
    \end{split}
    \end{align}
    where the absolute value of the determinant of $A$ restricted to $V^{i}$ is given by the product of the singular values of $A$ restricted to the vector subspace $V^{i}.$
\end{definition}

\begin{remark}
    Consider $A\in \text{GL}(\mathbb{E}_x,\mathbb{E}_y)$ and $A^{*}$ its adjoint as defined in (\ref{adjoint matrix}). Let $v_{1},\dots, v_{d}$ be an orthonormal basis of eigenvectors for $A^{*}A$ with corresponding eigenvalues $\sigma_{1}(A)^{2}\geq\dots\geq \sigma_{2}(A)^{2}$. For $i\in\{1,\dots, d\}$, define the vector subspace $S^{i}=\mathrm{span}\{v_{1},\dots, v_{i}\}$ and call $\mathcal{S}(A)=(S^{i})$ a \textit{singular flag} of $A$. 
    
    Let $f$ be a linear Schur-convex functional. We see that
    \begin{align}
    \begin{split}
    f(\kappa^{g}(A))&=\sum_{i=1}^{d}a_i\log\norm{\mathsf{\Lambda}^{i}A}\\
    &=\sum_{i=1}^{d}a_i\log|\det(A|S^{i})|=\tilde{f}(A, \mathcal{S}),
    \end{split}
    \end{align}
    Furthermore, we also have the equality
    \begin{equation}
    f(\kappa^{g}(A))=\max_{\mathcal{V}\in\mathfrak{F}_{x}}\tilde{f}(A,\mathcal{V}).
    \end{equation}
\end{remark}

Since we need to understand the evolution of the bundle map when restricted to a vector subspace, we define the following sequence using the induced functions on complete flags:  
\begin{equation} \label{seqence restricted to flags}
\tilde{F}_{n}(x,\mathcal{V})=\tilde{f}(\Phi_x^{n}, \mathcal{V})=\sum_{i=1}^{d}a_i\log|\det(\Phi^{n}_x|V_i)|,
\end{equation}
where $\mathcal{V}\in \mathfrak{F}_x$ is a complete flag on $\mathbb{E}_{x}.$ 

Taking $u_1,\dots, u_i$ an orthonormal basis for $V^i$, we get the following description of (\ref{seqence restricted to flags}):  
\begin{align} \label{determinants}
    \begin{split}
    \tilde{F}_{n}(x,\mathcal{F})&=\sum_{i=1}^{d}a_i\log\norm{\mathsf{\Lambda}^{i}\Phi^{n}_{x}(u_1\wedge\dots\wedge u_i)}\\
    &=\sum_{i=1}^{d}a_i\sum_{j=0}^{n-1}\log\dfrac{\norm{\mathsf{\Lambda}^{i}\Phi_{T^jx}(\Phi^{j}_xu_1\wedge\dots\wedge \Phi^{j}_xu_i)}}{\norm{\Phi^{j}_xu_1\wedge\dots\wedge \Phi^{j}_xu_i}}\\
    &=\sum_{i=1}^{d}a_i\sum_{j=0}^{n-1}\log |\det(\Phi_{T^jx}|\Phi^{j}_xV_i)|.
    \end{split}
\end{align}

\subsection{Evaluation of Schur-convex functionals at iterations of the bundle map at nearby points.}

Recall that $T:X\to X$ is a hyperbolic mixing homeomorphism, $0<\theta\leq 1$ is the Hölder exponent related to the bundle map $\Phi$ which covers $T$ and $\tau>0$ is the constant of expansion. 

Since it will be necessary to compare the growth of the fiber map in different fibers, we recall the following result from Kalinin and Sadovskaya \cite[Section~2.2]{KalininSadovskaya2013} will allow us to move objects between fibers. 

\begin{lemma} \label{transportmaps}
    Let $\mathbb{E}$ be a $\theta$-H\"older vector bundle endowed with a Riemannian norm $g$. Then there exists $\rho>0$ so that for all elements $x,y\in X$ where  $d(x,y)<\rho$, there are linear maps $I_{y\xleftarrow{}x}:\mathbb{E}_x\to\mathbb{E}_y$ with the following properties: 
    \begin{enumerate}
        \item If $x=y$, then $I_{x\xleftarrow{}x}$ is equal to the identity. 
        \item $I_{x\xleftarrow{}y}^{-1}=I_{y\xleftarrow{}x}$
        \item The map $I_{y\xleftarrow{}x}$ is $\theta$-H\"older on $x$ and $y$ for nearby points. 
    \end{enumerate}
\end{lemma}
Furthermore, these transition maps give a characterization of $\theta-$Hölder bundle maps.

\begin{lemma} \label{LemmaTransitionHolder}
    A bundle map $\Phi:\mathbb{E}\to\mathbb{E}$ covering a Lipschitz homeomorphism $T:X\to X$ is $\theta-$Hölder if and only if there is $C>0$ such that
    \begin{equation}
    \norm{I_{Ty\xleftarrow{}Tx}\circ\Phi_{x}-\Phi_{y}\circ I_{y\xleftarrow{}x}}\leq Cd(x,y)^{\theta}
    \end{equation}
    for all nearby $x,y\in X$.
\end{lemma}

Before we can enunciate the main lemma proved in this subsection, we will require the following elements. Fix $x, y\in X$, $k\in\mathbb{N}$ and $0<\delta<\rho$ so that for $i\in \{1,\dots, k\}$
\begin{equation}
    d(T^{i}x, T^{i}y)<\delta e^{-\tau(k-i)}. 
\end{equation}

Additionally, define the \textit{condition number} at $x\in X$ as
\begin{equation}
\text{cn}(\Phi_{x})=\norm{\Phi_{x}}\norm{\Phi_{x}^{-1}}.
\end{equation}
Furthermore, fix $\epsilon>0$ and define 
\begin{equation} \label{definitionB}
L(\epsilon)=\max\left\{\max_{x\in X}\text{cn}(\Phi_{x}),e^{\tau\theta+\epsilon} \right\}. 
\end{equation}

The following lemma will be the main technical result used in the proof of Theorem \ref{Reformulation of main theorem}. 

\begin{lemma} \label{difference abstract}
Assume that $T:X\to X$ is a hyperbolic mixing homeomorphism and $\Phi:\mathbb{E}\to\mathbb{E}$ is a $\theta$-Hölder bundle map covering $T$. Consider a nonzero vector space $W\subseteq \mathbb{E}_y$. Then for every $\epsilon>0$ there exist $C_{\epsilon}>0$ such that  
\begin{align}
    \left|\log\dfrac{|\det(I_{T^{k}x\xleftarrow[]{}T^{k}y}\circ \Phi^{k}_{x}\circ I_{x\xleftarrow[]{}y}|W)|}{|\det(\Phi^{k}_{y}|W)|}\right|
    \leq C_{\epsilon}\delta^{\theta}e^{k(\log{L(\epsilon)}-\tau\theta)}. 
\end{align}
\end{lemma}
In order to prove this result, we will need some preliminary tools and lemmas, which we will explain next. 

For $0\leq i\leq  k$, we will consider the invertible matrices 
\begin{equation}
A_i=\Phi_{T^{i}y},
\end{equation}
\begin{equation}
B_i=I_{T^{i+1}y\xleftarrow{}T^{i+1}x}\circ\Phi_{T^{i}x}\circ I_{T^{i}x\xleftarrow{}T^{i}y}.
\end{equation}

Note that by the condition that $0<\delta<\rho$, we can ensure the existence of $I_{T^{i}x\xleftarrow{}T^{i}y}$ for all $0\leq i\leq k$. 

Due to Lemma \ref{LemmaTransitionHolder}, we see that 
\begin{align} \label{norm between matrices}
\begin{split}
  \norm{A_{i}-B_{i}}&\leq C d(T^{i}x,T^{i}y)^{\theta}\\
  &\leq C \delta^{\theta}e^{-\tau\theta(k-i)}
  \end{split}
  \end{align}
  where $C>0$ comes from the H\"older continuity.

Fix a unit vector $v_{0}\in \mathbb{E}_{q_0}$ and set $w_0=v_0$. For $i\in \{0,\dots, k\}$ we define
    \begin{align}
    \begin{split}
        v_i&=\dfrac{A_{i-1}\dots A_0v_0}{\norm{A_{i-1}\dots A_0v_0}},\\
        w_i&=\dfrac{B_{i-1}\dots B_0w_0}{\norm{B_{i-1}\dots B_0w_0}}.
    \end{split}
    \end{align}

Furthermore, recall the gap metric $d_{\mathrm{gap}}:\mathbb{P(E}_{z})\times \mathbb{P(E}_{z})\to[0,\infty)$ defined on each $z\in X$ by the sine of the angle between the elements $\overline{u_1},\overline{u_2}\in \mathbb{P(E}_{z})$. That is
\begin{equation} \label{gap metric}
d_{\mathrm{gap}}(\overline{u_1},\overline{u_2})=\dfrac{\norm{u_1\wedge u_{2}}_{z}}{\norm{u_1}_x\norm{u_1}_z}.
\end{equation}
We see that for every $z_1,z_2\in X $,  $A,B\in \mathrm{GL}(\mathbb{E}_{z_1},\mathbb{E}_{z_2} )$ and $v,w\in \mathbb{E}_{z_1}$, we obtain the inequality 
    \begin{align} \label{bounds metric projective}
    \begin{split}
        d_{\mathrm{gap}}(\overline{A v},\overline{Bv})&=\dfrac{\norm{Av\wedge Bv}}{\norm{Av}\norm{Bv}}\\
        &=\dfrac{\norm{[Bv-(B-A)v]\wedge Bv}}{\norm{Av}\norm{Bv}}\\
        &\leq\dfrac{\norm{B-A}}{\norm{Av}}
    \end{split}
    \end{align}

Additionally, from \cite[Proposition~3.1]{BochiGaribaldi}, we know that 
\begin{equation} \label{Lipschitz condition number}
    d_{\mathrm{gap}}(Av,Aw)\leq \mathrm{cn}(A)d(v,w). 
\end{equation}

Hence, we obtain the following lemma. 

\begin{lemma} \label{DistanceVectorsInduction}
 Consider $x,y\in X$, $v_{i},w_{i}\in \mathbb{E}_{T^{i}y}$ and $A_i, B_i\in \mathrm{GL}(\mathbb{E}_{T^{i}y},\mathbb{E}_{T^{i+1}y})$ as defined above for $i\in\{1,\dots, k\}$. Then, for $m\in \{0,\dots, k\}$
\begin{align}
d_{\mathrm{gap}}(\overline{v_m},\overline{w_m})\leq \sum_{i=0}^{m-1}\dfrac{\norm{A_i-B_i}}{\norm{A_iw_i}}\prod_{j=i+1}^{m-1}\mathrm{cn}(A_j),
\end{align}
where we set that $\prod_{i=s}^{l}\mathrm{cn}(A_i)=1$ if $l<s$. 
\end{lemma}
\begin{proof}
    Let us show it by induction for the case $m\geq 0$. It is obvious for $d_{\mathrm{gap}}(\overline{v_0},\overline{w_0}).$
    Then, for $1\leq m\leq k-1$, using the inequalities (\ref{bounds metric projective}) and (\ref{Lipschitz condition number})
     \begin{align}
     \begin{split}
        \hspace{-0.2cm}d_{\mathrm{gap}}(\overline{v_{m+1}},\overline{w_{m+1}})&\leq d_{\mathrm{gap}}(\overline{A_mv_{m}}, \overline{A_mw_m})+d_{\mathrm{gap}}(\overline{A_mw_m},\overline{B_m w_{m}})\\
        &\hspace{-0.2cm}\leq \text{cn}(A_m)d(\overline{v_m},\overline{w_m})+\dfrac{\norm{B_m-A_m}}{\norm{A_mw_m}}\\
        &\hspace{-0.2cm}\leq \text{cn}(A_m)\left(\sum_{i=0}^{m-1}\dfrac{\norm{A_i-B_i}}{\norm{A_iw_i}}\prod_{j=i+1}^{m-1}\text{cn}(A_j)\right)+\dfrac{\norm{B_m-A_k}}{\norm{A_mw_m}}\\
        &\hspace{-0.2cm}= \sum_{i=0}^{m-1}\dfrac{\norm{A_i-B_i}}{\norm{A_iw_i}}\prod_{j=i+1}^{m}\text{cn}(A_j)+\dfrac{\norm{B_m-A_m}}{\norm{A_mw_m}}\\
        &\hspace{-0.2cm}= \sum_{i=0}^{m}\dfrac{\norm{A_i-B_i}}{\norm{A_iw_i}}\prod_{j=i+1}^{m}\text{cn}(A_j).
        \end{split}
    \end{align}
\end{proof}

To simplify the notation, for functions $\phi_1, \phi_2:\{0,\dots, k\}\to\mathbb{R}$, we write that $\phi_1(m)\preceq \phi_2(m)$ if there exist a positive constant $C=C(T, \theta, \epsilon)$ such that $\phi_1(m)\leq C\phi_2(m)$ for all $m \in \{0,\dots, k\}$. 

\begin{lemma} \label{lemmavectors}
For any $0\leq m\leq k$ we obtain
\begin{align}
\begin{split}
d_{\mathrm{gap}}(\overline{v_m},\overline{w_m})\preceq \exp\big(\theta\log{\delta}-\theta\tau k+m\log{L}\big).
\end{split}
\end{align}
\end{lemma}
\begin{proof}
    Since $X$ is a compact set and $\Phi$ is a $\theta-$H\"older bundle map, we have that for any $z\in X$ and $u\in \mathbb{E}_z$, $\norm{\Phi_z u}$ is bounded from below. Therefore, from Lemma \ref{DistanceVectorsInduction}, inequality (\ref{norm between matrices}) and the definition of $L$, we obtain for $m\in\{0,\dots, d-1\}$
    \begin{align}
    \begin{split}
        d_{\mathrm{gap}}(\overline{v_{m+1}},\overline{w_{m+1}})&\preceq \sum_{i=0}^{m}\exp\big(\theta\log{\delta}-\tau\theta(k-i)+(m-i)\log{L}\big)\\
        &= \exp\big(\theta\log{\delta}-\theta\tau k+m\log{L}\big)\sum_{i=0}^{m}\exp[i(\theta\tau-\log{L})]\\
        &\leq \exp\big(\theta\log{\delta}-\theta\tau k+m\log{L}\big)\sum_{i=0}^{\infty}\exp[i(\theta\tau-\log{L})].
    \end{split}
    \end{align}    
      Since $L\geq e^{\theta\tau+\epsilon}$
      \begin{align}
      \begin{split}
        d_{\mathrm{gap}}(\overline{v_{m+1}},\overline{w_{m+1}})&\leq\exp\big(\theta\log{\delta}-\theta\tau k+m\log{L}\big)\frac{1}{1-\exp(\theta\tau-\log{L})}\\
        &\preceq \exp\big(\theta\log{\delta}-\theta\tau k+m\log{L}\big).
        \end{split}
        \end{align}\qedhere
\end{proof}
The gap metric defined in (\ref{gap metric}) can be generalized to higher dimensional spaces. Indeed, as it is shown in \cite[Appendix~A]{anosovrepresentations}, for $z\in X$ and vector subspaces $U_{1},U_{2}\subseteq \mathbb{E}_{z}$ of the same dimension, the gap metric is given by
\begin{equation}
d_{\mathrm{gap}}(U_1,U_2)=\max\left\{\max_{a\in U_1^{\times}} \min_{b\in U_{2}^{\times}}\sin\angle(a,b),\max_{b\in U_2^{\times}} \min_{a\in U_{1}^{\times}}\sin\angle(a,b)\right\},
\end{equation}
where the superindex $\times$ means that we consider the vector spaces without the origin. Moreover, this formula can be simplified to

\begin{equation}
d_{\mathrm{gap}}(U_1,U_2)=\max_{a\in U_1^{\times}}\sin\angle(a,U_2).
\end{equation}
We can now proceed to prove the main result of this subsection

\begin{proof}[Proof of Lemma \ref{difference abstract}]
    From (\ref{determinants}), we see that 
\begin{align}
    \begin{split}
      &\left|\log|\det(I_{T^{k}x\xleftarrow[]{}T^{k}y}\circ \Phi^{k}_{x}\circ I_{x\xleftarrow[]{}y}|W)|-\log |\det(\Phi^{k}_{y}|W)|\right|\\
      &\leq \sum_{m=0}^{k}|\log|\det(B_m|B_{m-1}\dots B_{0}W)|-\log|\det(A_m|A_{m-1}\dots A_{0}W)||\\
      &\leq \sum_{m=0}^{k}|\log|\det(B_m|B_{m-1}\dots B_{0}W)|-\log|\det(A_m|B_{m-1}\dots B_{0}W)||\\
      &\hspace{0.5cm} +\sum_{m=0}^{k}|\log|\det(A_m|B_{m-1}\dots B_{0}W)|-\log|\det(A_m|A_{m-1}\dots A_{0}W)||.
      \end{split}
     \intertext{Since log is a Lipchitz function}
     \begin{split}
     &\left|\log|\det(I_{T^{k}x\xleftarrow[]{}T^{k}y}\circ \Phi^{k}_{x}\circ I_{x\xleftarrow[]{}y}|W)|-\log |\det(\Phi^{k}_{y}|W)|\right|\\
      &\preceq \sum_{m=0}^{k}||\det(B_m|B_{m-1}\dots B_{0}W)|-|\det(A_m|B_{m-1}\dots B_{0}W)||\\
      &\hspace{0.5cm} +\sum_{m=0}^{k}||\det(A_m|B_{m-1}\dots B_{0}W)|-|\det(A_m|A_{m-1}\dots A_{0}W)||.      
      \end{split}
      \end{align}
      Since the singular values are Lipchitz functions both on their dependence on the linear operator and on the linear subspace
      
      \begin{align}
      \begin{split}
      &\left|\log|\det(I_{T^{k}x\xleftarrow[]{}T^{k}y}\circ \Phi^{k}_{x}\circ I_{x\xleftarrow[]{}y}|W)|-\log |\det(\Phi^{k}_{y}|W)|\right|\\
      &\preceq \sum_{m=0}^{k}\norm{B_m-A_m}+\sum_{m=0}^{k}d_{\mathrm{gap}}(B_{m-1}\dots B_{0}W, A_{m-1}\dots A_{0}W)\\
      &=\sum_{m=0}^{k}\norm{B_m-A_m}+\max_{w\in W}\sin\angle(B_{m-1}\dots B_{0}w, A_{m-1}\dots A_{0}W)\\
      &\leq \sum_{m=0}^{k} \norm{B_m-A_m}+\max_{w\in W}\sin\angle(B_{m-1}\dots B_{0}w, A_{m-1}\dots A_{0}w). 
      \end{split}
      \intertext{Then, from (\ref{norm between matrices}) and Lemma \ref{lemmavectors} applied to the 
      first and second term respectively, we obtain}
      \begin{split}
      &\left|\log|\det(I_{T^{k}x\xleftarrow[]{}T^{k}y}\circ \Phi^{k}_{x}\circ I_{x\xleftarrow[]{}y}|W)|-\log |\det(\Phi^{k}_{y}|W)|\right|\\
      &\preceq \sum_{m=0}^{k} \exp\big(\theta\log{\delta}-\tau\theta(k-m)\big)+\exp\big(\theta\log{\delta}-\theta\tau k+m\log{L}\big)\\
      &= \exp\big(\theta\log{\delta}-k\tau\theta )\sum_{m=0}^{k} \exp\big(m\tau\theta \big)+\exp\big(m\log{L}\big)\\
      &\preceq \exp\big(\theta\log{\delta}-k\tau\theta )\sum_{m=0}^{k}\exp\big(m\log{L}\big)\\
      &\preceq \exp\big(\theta\log{\delta}-k\tau\theta )\left[\dfrac{\exp((k+1)\log{L})-1}{\exp(\log{L})-1}\right]\\
      &\preceq \delta^{\theta}\exp\big[k(\log{L}-\tau\theta)\big]. 
      \end{split}
  \end{align}  
  Since the constant that is disregarded by the inequality symbol $\preceq$ can only depend on $T, \theta$ and $\epsilon$, the result follows. 
\end{proof}

\begin{remark} \label{process also works backwards}
    Note that the same result can be obtained with the fiber map $\Phi^{-1}:\mathbb{E}\to\mathbb{E}$ over the hyperbolic homeomorphism $T^{-1}:X\to X$. 
\end{remark}

\subsection{Approximating values of the functional using periodic points}

In order to prove Theorem \ref{Reformulation of main theorem}, we will show that for every $x\in X$ we can approximate the value of $F_{n}(x)$ with the value of $F_{n}$ evaluated at a periodic point. Moreover, from (\ref{F in terms of wedge}), we see that we only have to study the exterior products of the fiber map $\norm{\mathsf{\Lambda}^{i}\Phi^{n}}$. 

Fix $x_{0}\in X$. Take $n\gg 1$ large, $\epsilon>0$ as in (\ref{definitionB})  and $\gamma\in\mathbb{Q}$ such that 
\begin{equation}
    \gamma\geq\max\left\{\max_{x\in X}\dfrac{\log{\text{cn}}(\Phi_{x})}{\tau\theta},1 \right\}+\frac{\epsilon}{\tau\theta}.
\end{equation}
Due to Proposition \ref{change of uniform spectrum due to iterations}, our results will not be affected by imposing the condition that $n\gamma/2 \in \mathbb{N}$.

Consider $\tilde{x}=T^{-\gamma n/2+n/2}x_0$. Then by Lemma \ref{approximation} there exists a periodic point $\tilde{p}\in X$ such that $T^{\gamma n}\tilde{p}=\tilde{p}$ and which exponentially shadows the orbit segment  $\tilde{x}, \dots, T^{\gamma n}\tilde{x}$. Hence, there is $C>0$ such that for $i\in \{0,\dots, n\gamma\}$
\begin{equation}
    d(T^{i}\tilde{x}, T^{i}\tilde{p})\leq Ce^{-\tau\min\{i,n\gamma-i\}}. 
\end{equation}

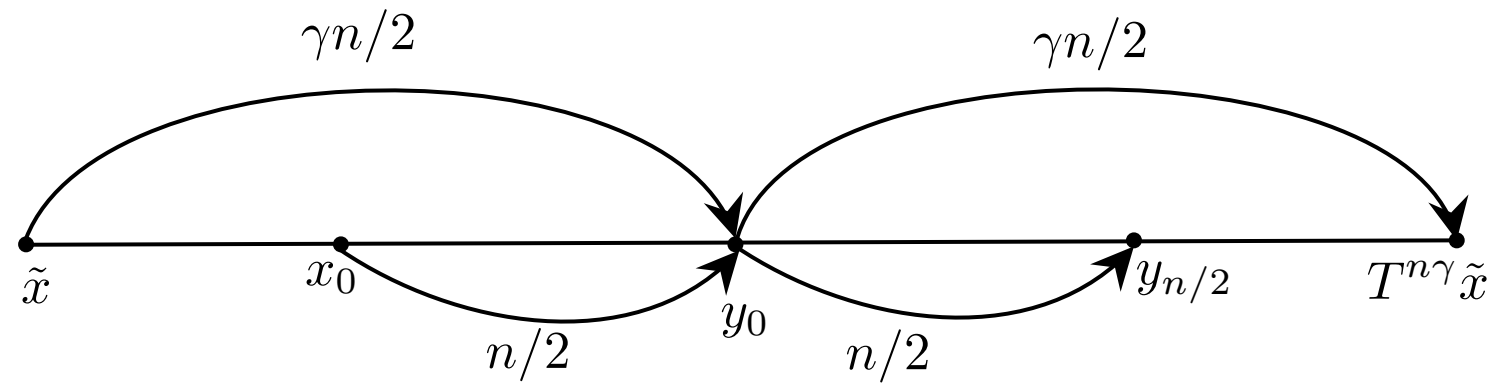
\begin{figure}[h] 

\centering

\tikzset{every picture/.style={line width=0.75pt}} 

\begin{tikzpicture}[x=0.75pt,y=0.75pt,yscale=-1,xscale=1]

\draw    (301.4,100.6) -- (660.9,99.5) ;
\draw  [fill={rgb, 255:red, 0; green, 0; blue, 0 }  ,fill opacity=1 ] (300.4,100.5) .. controls (300.4,99.67) and (301.07,99) .. (301.9,99) .. controls (302.73,99) and (303.4,99.67) .. (303.4,100.5) .. controls (303.4,101.33) and (302.73,102) .. (301.9,102) .. controls (301.07,102) and (300.4,101.33) .. (300.4,100.5) -- cycle ;
\draw  [fill={rgb, 255:red, 0; green, 0; blue, 0 }  ,fill opacity=1 ] (379.4,100.5) .. controls (379.4,99.67) and (380.07,99) .. (380.9,99) .. controls (381.73,99) and (382.4,99.67) .. (382.4,100.5) .. controls (382.4,101.33) and (381.73,102) .. (380.9,102) .. controls (380.07,102) and (379.4,101.33) .. (379.4,100.5) -- cycle ;
\draw  [fill={rgb, 255:red, 0; green, 0; blue, 0 }  ,fill opacity=1 ] (478.4,100.5) .. controls (478.4,99.67) and (479.07,99) .. (479.9,99) .. controls (480.73,99) and (481.4,99.67) .. (481.4,100.5) .. controls (481.4,101.33) and (480.73,102) .. (479.9,102) .. controls (479.07,102) and (478.4,101.33) .. (478.4,100.5) -- cycle ;
\draw  [fill={rgb, 255:red, 0; green, 0; blue, 0 }  ,fill opacity=1 ] (578.4,99.5) .. controls (578.4,98.67) and (579.07,98) .. (579.9,98) .. controls (580.73,98) and (581.4,98.67) .. (581.4,99.5) .. controls (581.4,100.33) and (580.73,101) .. (579.9,101) .. controls (579.07,101) and (578.4,100.33) .. (578.4,99.5) -- cycle ;
\draw  [fill={rgb, 255:red, 0; green, 0; blue, 0 }  ,fill opacity=1 ] (659.4,99.5) .. controls (659.4,98.67) and (660.07,98) .. (660.9,98) .. controls (661.73,98) and (662.4,98.67) .. (662.4,99.5) .. controls (662.4,100.33) and (661.73,101) .. (660.9,101) .. controls (660.07,101) and (659.4,100.33) .. (659.4,99.5) -- cycle ;
\draw    (301.9,99) .. controls (320.71,50.54) and (461.03,49.2) .. (479.16,96.79) ;
\draw [shift={(479.9,99)}, rotate = 253.71] [fill={rgb, 255:red, 0; green, 0; blue, 0 }  ][line width=0.08]  [draw opacity=0] (10.72,-5.15) -- (0,0) -- (10.72,5.15) -- (7.12,0) -- cycle    ;
\draw    (479.9,100.5) .. controls (492.48,48.62) and (646.28,49.79) .. (660.36,97.3) ;
\draw [shift={(660.9,99.5)}, rotate = 258.73] [fill={rgb, 255:red, 0; green, 0; blue, 0 }  ][line width=0.08]  [draw opacity=0] (10.72,-5.15) -- (0,0) -- (10.72,5.15) -- (7.12,0) -- cycle    ;
\draw    (380.9,102) .. controls (409.27,121.18) and (453.71,129.26) .. (478.99,104.01) ;
\draw [shift={(480.9,102)}, rotate = 131.99] [fill={rgb, 255:red, 0; green, 0; blue, 0 }  ][line width=0.08]  [draw opacity=0] (10.72,-5.15) -- (0,0) -- (10.72,5.15) -- (7.12,0) -- cycle    ;
\draw    (479.9,101) .. controls (508.27,120.18) and (552.71,128.26) .. (577.99,103.01) ;
\draw [shift={(579.9,101)}, rotate = 131.99] [fill={rgb, 255:red, 0; green, 0; blue, 0 }  ][line width=0.08]  [draw opacity=0] (10.72,-5.15) -- (0,0) -- (10.72,5.15) -- (7.12,0) -- cycle    ;

\draw (299.13,104.8) node [anchor=north west][inner sep=0.75pt]   [align=left] {$\displaystyle \tilde{x}$};
\draw (370.4,103.5) node [anchor=north west][inner sep=0.75pt]   [align=left] {$\displaystyle x_{0}$};
\draw (474.6,114) node [anchor=north west][inner sep=0.75pt]   [align=left] {$\displaystyle y_0$};
\draw (578.8,103.3) node [anchor=north west][inner sep=0.75pt]   [align=left] {$\displaystyle y_{n/2}$};
\draw (636.8,103.13) node [anchor=north west][inner sep=0.75pt]   [align=left] {$\displaystyle T^{n\gamma } \tilde{x}$};
\draw (369.33,40) node [anchor=north west][inner sep=0.75pt]   [align=left] {$\displaystyle \gamma n/2$};
\draw (415.67,120) node [anchor=north west][inner sep=0.75pt]   [align=left] {$\displaystyle n/2$};
\draw (553,42) node [anchor=north west][inner sep=0.75pt]   [align=left] {$\displaystyle \gamma n/2$};
\draw (506,120.33) node [anchor=north west][inner sep=0.75pt]   [align=left] {$\displaystyle n/2$};

\end{tikzpicture}
\caption{Orbit of $x_0\in X$.}
\label{labeledOrbit}
\end{figure}

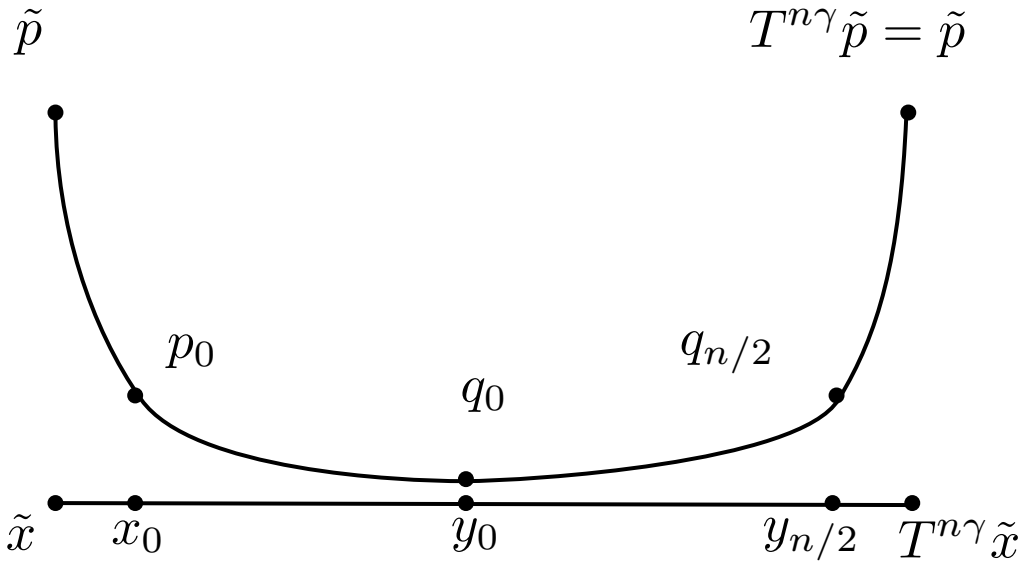
\begin{figure}[h] 

\centering

\tikzset{every picture/.style={line width=0.75pt}} 

\begin{tikzpicture}[x=0.75pt,y=0.75pt,yscale=-1,xscale=1]

\draw    (366.4,187.5) -- (582.9,188) ;
\draw  [fill={rgb, 255:red, 0; green, 0; blue, 0 }  ,fill opacity=1 ] (366.4,89.5) .. controls (366.4,88.67) and (367.07,88) .. (367.9,88) .. controls (368.73,88) and (369.4,88.67) .. (369.4,89.5) .. controls (369.4,90.33) and (368.73,91) .. (367.9,91) .. controls (367.07,91) and (366.4,90.33) .. (366.4,89.5) -- cycle ;
\draw  [fill={rgb, 255:red, 0; green, 0; blue, 0 }  ,fill opacity=1 ] (386.4,160.5) .. controls (386.4,159.67) and (387.07,159) .. (387.9,159) .. controls (388.73,159) and (389.4,159.67) .. (389.4,160.5) .. controls (389.4,161.33) and (388.73,162) .. (387.9,162) .. controls (387.07,162) and (386.4,161.33) .. (386.4,160.5) -- cycle ;
\draw  [fill={rgb, 255:red, 0; green, 0; blue, 0 }  ,fill opacity=1 ] (469.4,181.5) .. controls (469.4,180.67) and (470.07,180) .. (470.9,180) .. controls (471.73,180) and (472.4,180.67) .. (472.4,181.5) .. controls (472.4,182.33) and (471.73,183) .. (470.9,183) .. controls (470.07,183) and (469.4,182.33) .. (469.4,181.5) -- cycle ;
\draw  [fill={rgb, 255:red, 0; green, 0; blue, 0 }  ,fill opacity=1 ] (562.4,160.5) .. controls (562.4,159.67) and (563.07,159) .. (563.9,159) .. controls (564.73,159) and (565.4,159.67) .. (565.4,160.5) .. controls (565.4,161.33) and (564.73,162) .. (563.9,162) .. controls (563.07,162) and (562.4,161.33) .. (562.4,160.5) -- cycle ;
\draw  [fill={rgb, 255:red, 0; green, 0; blue, 0 }  ,fill opacity=1 ] (580.4,89.5) .. controls (580.4,88.67) and (581.07,88) .. (581.9,88) .. controls (582.73,88) and (583.4,88.67) .. (583.4,89.5) .. controls (583.4,90.33) and (582.73,91) .. (581.9,91) .. controls (581.07,91) and (580.4,90.33) .. (580.4,89.5) -- cycle ;
\draw  [fill={rgb, 255:red, 0; green, 0; blue, 0 }  ,fill opacity=1 ] (366.4,187.5) .. controls (366.4,186.67) and (367.07,186) .. (367.9,186) .. controls (368.73,186) and (369.4,186.67) .. (369.4,187.5) .. controls (369.4,188.33) and (368.73,189) .. (367.9,189) .. controls (367.07,189) and (366.4,188.33) .. (366.4,187.5) -- cycle ;
\draw  [fill={rgb, 255:red, 0; green, 0; blue, 0 }  ,fill opacity=1 ] (386.4,187.5) .. controls (386.4,186.67) and (387.07,186) .. (387.9,186) .. controls (388.73,186) and (389.4,186.67) .. (389.4,187.5) .. controls (389.4,188.33) and (388.73,189) .. (387.9,189) .. controls (387.07,189) and (386.4,188.33) .. (386.4,187.5) -- cycle ;
\draw  [fill={rgb, 255:red, 0; green, 0; blue, 0 }  ,fill opacity=1 ] (469.4,187.5) .. controls (469.4,186.67) and (470.07,186) .. (470.9,186) .. controls (471.73,186) and (472.4,186.67) .. (472.4,187.5) .. controls (472.4,188.33) and (471.73,189) .. (470.9,189) .. controls (470.07,189) and (469.4,188.33) .. (469.4,187.5) -- cycle ;
\draw  [fill={rgb, 255:red, 0; green, 0; blue, 0 }  ,fill opacity=1 ] (561.4,187.5) .. controls (561.4,186.67) and (562.07,186) .. (562.9,186) .. controls (563.73,186) and (564.4,186.67) .. (564.4,187.5) .. controls (564.4,188.33) and (563.73,189) .. (562.9,189) .. controls (562.07,189) and (561.4,188.33) .. (561.4,187.5) -- cycle ;
\draw  [fill={rgb, 255:red, 0; green, 0; blue, 0 }  ,fill opacity=1 ] (581.4,187.5) .. controls (581.4,186.67) and (582.07,186) .. (582.9,186) .. controls (583.73,186) and (584.4,186.67) .. (584.4,187.5) .. controls (584.4,188.33) and (583.73,189) .. (582.9,189) .. controls (582.07,189) and (581.4,188.33) .. (581.4,187.5) -- cycle ;
\draw    (367.9,90.5) .. controls (368.2,117.6) and (376.6,143.4) .. (389.4,161.5) .. controls (402.03,179.36) and (453.8,182.02) .. (470.25,182) .. controls (470.47,182) and (470.69,182) .. (470.9,182) .. controls (486.59,181.91) and (554.6,177.6) .. (564.4,161.5) .. controls (564.45,161.42) and (564.5,161.34) .. (564.54,161.26) .. controls (564.56,161.23) and (564.58,161.2) .. (564.6,161.17) .. controls (574.28,145.14) and (579.81,126.55) .. (581.4,90.5) ;

\draw (356.33,62) node [anchor=north west][inner sep=0.75pt]   [align=left] {$\displaystyle \tilde{p}$};
\draw (540.4,62) node [anchor=north west][inner sep=0.75pt]   [align=left] {$\displaystyle T^{n\gamma } \tilde{p}=\tilde{p}$};
\draw (394.4,143.8) node [anchor=north west][inner sep=0.75pt]   [align=left] {$\displaystyle p_{0}$};
\draw (468,154.8) node [anchor=north west][inner sep=0.75pt]   [align=left] {$\displaystyle q_0$};
\draw (523,143) node [anchor=north west][inner sep=0.75pt]   [align=left] {$\displaystyle  q_{n/2}$};
\draw (353.93,188.8) node [anchor=north west][inner sep=0.75pt]   [align=left] {$\displaystyle \tilde{x}$};
\draw (380.4,190.5) node [anchor=north west][inner sep=0.75pt]   [align=left] {$\displaystyle x_{0}$};
\draw (465.6,189.8) node [anchor=north west][inner sep=0.75pt]   [align=left] {$\displaystyle y_0$};
\draw (543.8,189.8) node [anchor=north west][inner sep=0.75pt]   [align=left] {$\displaystyle y_{n/2}$};
\draw (577.8,190.13) node [anchor=north west][inner sep=0.75pt]   [align=left] {$\displaystyle T^{n\gamma } \tilde{x}$};

\end{tikzpicture}
\caption{Exponential approximation of the orbit of $x_0\in X$ by the periodic point $p_0$.}
\label{ExponentialApproximation}
\end{figure}

Denote $p_0=T^{\gamma n/2-n/2}\tilde{p}$. Write $q_0=T^{n/2}p_0$ and $y_{0}=T^{n/2}x_0$. Additionally, write $q_{i}=T^{i}q_0$ and $y_{i}=T^{i}y_0$ for $- n/2\leq i\leq n/2 $. See Figure \ref{labeledOrbit} for a representation of the orbit of $x_0$ and Figure \ref{ExponentialApproximation} for a portrayal of the behavior of $p_0.$

Consider $f:\mathfrak{a}_{+}^{d}\to \mathbb{R}$ a Schur-convex functional and $F_n$ as defined in  (\ref{Iteratedfunctionals}).  
We will proceed to find a lower bound for $\log\norm{\mathsf{\Lambda}^{j}\Phi^{(n)}_{x_0}}$ where $1\leq j\leq d$. To that end, consider a complete flag 
\begin{align}
\mathcal{V}_{-n/2}(p_0)=(V_{-n/2}^{j})\in \mathfrak{F}_{p_0}
\intertext{such that}
F_n(p_0)=\tilde{F}_n(p_0,\mathcal{V}_{-n/2}(p_0)), 
\end{align}
that is, $\mathcal{V}_{-n/2}(p_0)$ is a singular flag for $\Phi^{n}_{p_0}$. 

For each $j\in \{0,\dots, d\}$, write 
\begin{equation}
V_0^{j}=\Phi^{n/2}_{p_{0}}(V^{j}_{-n/2}),
\end{equation}
and consider the complete flag $\mathcal{V}_{0}(q_{0})=(V^{j}_{0})\in \mathfrak{F}_{q_{0}}$.

For each $j\in \{1,\dots, d\}$, define
\begin{equation}
W_{-n/2}^{j}=I_{p_{0}\xleftarrow[]{}x_{0}}\circ\Phi^{-n/2}_{y_0}\circ I_{y_0\xleftarrow[]{}q_0} (V^{j}_{0}). 
\end{equation}

Since we know that 
\begin{align}
\begin{split}
&\log\norm{\mathsf{\Lambda}^{j}I_{T^{n}p_{0}\xleftarrow[]{}T^{n}x_{0}}\circ\Phi^{n}_{x_0}\circ I_{x_{0}\xleftarrow[]{}p_{0}}}\\
&\hspace{0cm}\geq \log|\det(I_{T^{n}p_{0}\xleftarrow[]{}T^{n}x_{0}}\circ\Phi^{n}_{x_0}\circ I_{x_{0}\xleftarrow[]{}p_{0}}|W^{j}_{-n/2})|,
\end{split}
\end{align}
our strategy to find a lower bound for $\log\norm{\mathsf{\Lambda}^{j}\Phi^{n}(x_0)}$ will consist on showing that the difference 
\begin{equation} \label{division of substraction}
    |\log|\det(I_{T^{n}p_{0}\xleftarrow[]{}T^{n}x_{0}}\circ\Phi^{n}_{x_0}\circ I_{x_{0}\xleftarrow[]{}p_{0}}|W^{j}_{-n/2})|- \log|\det(\Phi^{n}(p_0)|V^{j}_{-n/2})||
\end{equation}
is bounded by a constant, so that we can find a lower estimate in terms of the periodic point.

If $t_1, \dots, t_j$ is an orthonormal basis of $V^{j}_{0}$, we see that
\begin{align}
\begin{split}
    &\log|\det(I_{T^{n}p_{0}\xleftarrow[]{}T^{n}x_{0}}\circ\Phi^{n}_{x_0}\circ I_{x_{0}\xleftarrow[]{}p_{0}}|W^{j}_{-n/2})|\\
    &=\log\left(\dfrac{\norm{\mathsf{\Lambda}^{j}I_{T^{n}p_{0}\xleftarrow[]{}T^{n}x_{0}}\circ\Phi^{n}_{x_0}\circ \Phi^{-n/2}_{y_0}\circ I_{y_{0}\xleftarrow[]{}q_{0}}(t_1\wedge\dots\wedge t_j)}}{\norm{\mathsf{\Lambda}^{j}I_{p_{0}\xleftarrow[]{}x_{0}}\circ\Phi^{-n/2}_{y_0}\circ I_{y_0\xleftarrow[]{}q_0}(t_1\wedge\dots\wedge t_j)}}\right)\\
    &=\log|\det(I_{T^{n}p_{0}\xleftarrow[]{}T^{n}x_{0}}\circ \Phi^{n/2}_{y_0}\circ I_{y_0\xleftarrow[]{}q_0}|V^{j}_{0})|\\
    &\hspace{2.5cm}-\log|\det(I_{p_{0}\xleftarrow[]{}x_{0}}\circ \Phi^{-n/2}_{y_0}\circ I_{y_0\xleftarrow[]{}q_0}|V^{j}_{0})|.
\end{split}
\end{align}
    We obtain a similar result for $p_0$
    \begin{equation}
    \log|\det(\Phi^{n}_{p_0}|V^{j}_{-n/2})|=\log|\det(\Phi^{n/2}_{q_0}|V^{j}_{0})|-\log|\det(\Phi^{-n/2}_{q_0}|V^{j}_{0})|.
    \end{equation}
Hence, we can split the difference (\ref{division of substraction}) as
\begin{align} \label{separate into substraction}
    \begin{split}
    &|\log|\det(I_{T^{n}p_{0}\xleftarrow[]{}T^{n}x_{0}}\circ\Phi^{n}_{x_0}\circ I_{x_{0}\xleftarrow[]{}p_{0}}|W^{j}_{-n/2})|-\log|\det(\Phi^{n}_{p_0}|V^{j}_{-n/2})||\\
    &\leq|\log|\det(I_{T^{n}p_{0}\xleftarrow[]{}T^{n}x_{0}}\circ \Phi^{n/2}_{y_0}\circ I_{y_0\xleftarrow[]{}q_0}|V^{j}_{0})|-\log|\det(\Phi^{n/2}_{q_0}|V^{j}_{0})||\\
    &\hspace{0.5cm}+|\log|\det(I_{p_{0}\xleftarrow[]{}x_{0}}\circ \Phi^{-n/2}_{y_0}\circ I_{y_0\xleftarrow[]{}q_0}|V^{j}_{0})|-\log|\det(\Phi^{-n/2}_{q_0}|V^{j}_{0})||.
    \end{split}
\end{align}

\begin{lemma} \label{difference}
For every $\epsilon>0$, there exist a positive constant $C=C(T, \theta, \epsilon )$ such that  
\begin{align}
    \left|\log\dfrac{|\det(I_{T^{n}p_{0}\xleftarrow[]{}T^{n}x_{0}}\circ\Phi^{n}_{x_0}\circ I_{x_{0}\xleftarrow[]{}p_{0}}|W^{j}_{-n/2})|}{|\det(\Phi^{n}_{p_0}|V^{j}_{-n/2})|}\right|
    \leq C.
\end{align}
\end{lemma}
\begin{proof}
    Since for all $i\in \{0,\dots, n/2\}$ the distances between the orbit segments of $y_{0}$ and $q_{0}$ are 
    \begin{equation}
        d(y_i,q_i)\leq Ce^{-\tau(n\gamma/2-i)}=Ce^{-\tau(n\gamma/2-n/2)}e^{-\tau(n/2-i)},
    \end{equation}
    we can apply Lemma \ref{difference abstract} with $\delta=e^{-\tau(n\gamma/2-n/2)}$ and $k=n/2$.
    
    Hence, since $\tau \theta \gamma\geq \log{L}$ we obtain that
    \begin{align}
        \begin{split}
        \left|\log\dfrac{|\det(I_{T^{n}p_{0}\xleftarrow[]{}T^{n}x_{0}}\circ \Phi^{n/2}_{y_0}\circ I_{y_0\xleftarrow[]{}q_0}|V^{j}_{0})|}{|\det(\Phi^{n/2}_{q_0}|V^{j}_{0})|}\right|
    & \leq C_{\epsilon}e^{-\tau \theta n\gamma/2+(n/2)\log{L}}\\
    & \leq C. 
    \end{split}
    \end{align}
    The lemma follows from doing the same process with the other term in the right hand side of inequality (\ref{separate into substraction}). 
\end{proof}

\subsection{Proof of Theorem \ref{Reformulation of main theorem}}

We will begin this section by showing that for any linear Schur-convex function $f$, $\max_{x\in X}f(\kappa^{g}(\Phi_{x}))$ is close to $\lambda_{f}$ under a change of Reimannian metric. 

\begin{prop} \label{NewMetric}
    Consider $f:\mathfrak{a}_{+}^{d}\to\mathbb{R}$ a linear Schur-convex functional and let $\epsilon>0$. Then there is a $C^{0}$ Riemannian metric $g^{\epsilon}:\mathbb{E}^{2}\to\mathbb{R}$ such that for all $x\in X$
    \begin{equation}
    f(\kappa^{g^{\epsilon}}(\Phi_{x}))\leq \lambda_{f}^{\top}+\epsilon,
    \end{equation}
    where $\kappa^{g^{\epsilon}}$ is the Cartan projection with respect to the metric $g^{\epsilon}$.
\end{prop}
\begin{proof}
    Take $\epsilon>0$ and $\mathcal{O}\subseteq \mathfrak{a}_{+}^{d}$ an open set such that $\mathrm{Maj}(\mathcal{U}(\Phi))\subseteq \mathcal{O}$ and 
    \begin{equation}
    \sup_{\mathcal{O}}f\leq \lambda_{f}^{\top}+\epsilon.
    \end{equation}
    Hence, by Proposition \ref{Change of norm}, there is a $C^{0}$ Riemannian metric $g^{\epsilon}$ such that
    $\kappa^{g^{\epsilon}}(\Phi)\subseteq \mathcal{O}$. Thus, for all $x\in X$
    \begin{equation}
     f(\kappa^{g^{\epsilon}}(\Phi_{x}))\leq \sup_{\mathcal{O}}f\leq \lambda_{f}^{\top}+\epsilon.
     \end{equation}
\end{proof}

Consider the Schur-convex functional $(\xi_1-\xi_d):\mathfrak{a}_{+}^{d}\to\mathbb{R}$ given by the vector $(1,0,\dots, 0, -1)\in \mathfrak{a}_{+}^{d}$. Note that $\log\text{cn}(\Phi_{x})=(\xi_1-\xi_d)(\kappa^{g}(\Phi_{x})).$ Thus, by Proposition \ref{NewMetric}, for any $\epsilon>0$ we can find a $C^{0}$ metric $g^{\epsilon}$ such that for any $x\in X$
\begin{equation}
(\xi_1-\xi_d)(\kappa^{g^{\epsilon}}(\Phi_x))\leq (\xi_1-\xi_d)^{\top}+\frac{\epsilon}{\theta\tau}. 
\end{equation}

Hence, by taking 
\begin{equation} 
    \gamma\geq\max\left\{\dfrac{(\xi_1-\xi_d)^{\top}}{\theta\tau}, 1 \right\}+\frac{\epsilon}{\theta\tau}, 
\end{equation}
we will be able to use Lemma \ref{difference} in the proof of Theorem \ref{Reformulation of main theorem}.

\begin{proof}[Proof of Theorem \ref{Reformulation of main theorem}]
Consider a sequence an increasing sequence of natural numbers $\{n_j\}\subseteq \mathbb{N}$ and a sequence $\{x_{j}\}\subset X$ such that 
\begin{equation}
\alpha_{f}^{\bot}=\lim_{j\to\infty}n_j^{-1}F_{n_j}(x_{j}).
\end{equation}
From Lemma $\ref{difference}$, we know that for each $x_{j}$ we can find a periodic point $p_{j}$ of period $n_j\gamma$ and complete flags $\mathcal{V}_{j}=(V_{n_j}^{j})_{j}\in\mathfrak{F}_{p_{j}}$ and $\mathcal{W}_{n_j}=(W_{n_j}^{j})_{j}\in\mathfrak{F}_{p_{j}}$,  where $\mathcal{V}_{n_j}$ is a singular flag so that $\tilde{F}_{n_j}(p_{j},\mathcal{V}_{n_{j}})=F_{n_j}(p_{j})$. Therefore
    \begin{align}
        \begin{split}
        F_{n_j}(x_{j})&=\sum_{i=1}^{d}a_i\log\norm{\mathsf{\Lambda}^{i}\Phi^{n_j}_{x_{j}}}\\
        &\geq\sum_{i=1}^{d}a_i\log\norm{\mathsf{\Lambda}^{i}I_{T^{n_j}p_{j}\xleftarrow[]{}T^{n_j}x_{j}}\circ\Phi^{n_j}_{x_{j}}\circ I_{x_{j}\xleftarrow[]{}p_{j}}}\\
        &\hspace{0.5cm}-\sum_{i=1}^{d}a_i\left(\log\norm{\mathsf{\Lambda}^{i}I_{T^{n_j}x_{j}\xleftarrow[]{}T^{n_j}p_{j}}}+\log\norm{\mathsf{\Lambda}^{i}I_{p_{j}\xleftarrow[]{}x_{j}}}\right).
        \end{split}
        \intertext{Since the transition maps $I_{z\xleftarrow[]{}y}$ are continuous, there exist a constant $C_1>0$ such that}
        \begin{split}
        F_{n_j}(x_{j})&\geq\sum_{i=1}^{d}a_i\log\norm{\mathsf{\Lambda}^{i}I_{T^{n_j}p_{j}\xleftarrow[]{}T^{n_j}x_{j}}\circ\Phi^{n_j}_{x_{j}}\circ I_{x_{j}\xleftarrow[]{}p_{j}}}-C_1\\
        &\hspace{-1cm}\geq\sum_{i=1}^{d}a_i\log|\det(I_{T^{n_j}p_{j}\xleftarrow[]{}T^{n_j}x_{j}}\circ \Phi^{n_j}_{x_{j}}\circ I_{x_{j}\xleftarrow[]{}p_{j}}|W^{i}_{p_{j}})|-C_1. 
        \end{split}
        \end{align}
        Hence, from Lemma \ref{difference}, there is $C_2>0$ such that
        \begin{align}
        \begin{split}
        F_{n_j}(x_{j})&\geq\sum_{i=1}^{d}a_i\log|\det(\Phi^{n_j}_{p_{j}}|V^{i}_{n_j})|-C_2\\
        &\geq \tilde{F}_{n_j}(p_{j},\mathcal{V}_{n_{j}})-C_2,
        \end{split}
        \end{align}
        and since $\mathcal{V}_{n_{j}}$ is a singular flag for $\Phi^{n_j}_{p_j}$
        \begin{equation}
        F_{n_j}(x_{j})\geq F_{n_j}(p_{j})-C_2.
        \end{equation}
        Thus, from (\ref{definition of top}), for $\delta>0$ and $n_j$ sufficiently large and recalling the definition of $F_{n}$ (\ref{Iteratedfunctionals})
        \begin{align}
        \begin{split}
        F_{n_j}(x_{j})&\geq F_{n_j\gamma}(p_{j})-(\lambda_{f}^{\top}+\delta)(n_j\gamma -n_j)-C_2\\
        &\geq f(\kappa^{g}(\Phi_{p_{j}}^{n_j\gamma}))-(\lambda_{f}^{\top}+\delta)(n_j\gamma -n_j)-C_2.
        \end{split}
        \intertext{Since $f$ is linear Schur-convex, due to Weyl's Theorem (\ref{EigenvaluesVSSingularValues})}
        F_{n_j}(x_{j})&\geq f(\zeta(\Phi_{p_{j}}^{n_j\gamma}))-(\lambda_{f}^{\top}+\delta)(n_j\gamma -n_j)-C_2.
        \intertext{Thus, multiplying by $n_j^{-1}$ and taking the limsup to both sides of the inequality}
        \begin{split}
        \hspace{-0.2cm}\lim_{j\to\infty}\dfrac{1}{n_j} F_{n_j}(x_{j})&\geq \limsup_{j\to\infty}\dfrac{1}{n_j}\left(f(\zeta(\Phi_{p_{j}}^{n_j\gamma}))- (\lambda_{f}^{\top}+\delta)(n_j\gamma -n_j)\right)\\
        &\geq \liminf_{j\to\infty}\dfrac{n_j\gamma}{n_j}f\left(\dfrac{1}{n_j\gamma}\zeta(\Phi_{p_{j}}^{n_j\gamma})\right)- (\lambda_{f}^{\top}+\delta)(\gamma -1).
        \end{split}
        \end{align}
        Then, by Proposition \ref{PeriodicPointsLyapunovSpectrum}, since $\frac{1}{n_j\gamma}\zeta(\Phi_{p_{n_j}}^{n_j\gamma})\in \mathcal{L}(\Phi)$
        \begin{equation}
        \alpha_{f}^{\bot}\geq \gamma\lambda_{f}^{\bot}- (\lambda_{f}^{\top}+\delta)(\gamma -1).
        \end{equation}
    Since this happens for all $\delta>0$, we conclude 
    \begin{equation}
    \lambda^{\top}_{f}-\alpha_{f}^{\bot}\leq \gamma(\lambda^{\top}_{f}-\lambda_{f}^{\bot}).
    \end{equation}
    \qedhere
\end{proof}
\section{Applications of the theorem} \label{Corollaries}
In the following section we will go over some theorems found in the literature which can be proved with the use of the uniform spectrum. 

The following statement is a quantitative version of \cite[Theorem~3.10]{DewittGogolev}. 

\begin{corollary} \label{CorollaryDeWittGogoglev}
    Assume that $T:X\to X$ is a hyperbolic mixing homeomorphism with constant of expansion $\tau>0$ and $\Phi:\mathbb{E}\to\mathbb{E}$ is a $\theta$-Hölder bundle map covering $T$. Suppose that $\mathcal{L}(\Phi)\subseteq \prod_{i=1}^{d}[\lambda_{i}-\delta,\lambda_{i}+\delta]$ where $\lambda=(\lambda_1,\dots, \lambda_d)\in \mathfrak{a}_{+}^{d}$ such that $\lambda_{m}>\lambda_{m+1}$ for some $1\leq m\leq d-1$ and $0<\delta< \frac{\lambda_m-\lambda_{m+1}}{(2^{m+1}+2^{m}-2)\gamma_{0}}$. Then $\Phi$ has a dominated splitting of index $m$. 
\end{corollary}

In order to obtain this result, we must find some bounds on $\mathcal{U}(\Phi)$ for $\mathcal{L}(\Phi)$ narrow. Hence, by Theorem \ref{MainTheorem}, it is enough to find bounds on $\mathrm{Aura}_{\gamma_{0}}(\mathcal{L}(\Phi))$.

\begin{lemma}
    Assume that $T:X\to X$ is a hyperbolic mixing homeomorphism with constant of expansion $\tau>0$ and $\Phi:\mathbb{E}\to\mathbb{E}$ is a $\theta$-Hölder bundle map covering $T$. Suppose that there exist $\delta>0$ such that  $\mathcal{L}(\Phi)\subseteq \prod_{i=0}^{d}[\lambda_ {i}-\delta, \lambda_{i}+\delta ]$ where $\lambda=(\lambda_1,\dots, \lambda_d)\in \mathfrak{a}_{+}^{d}$. Then, 
    \begin{equation}
            \mathrm{Aura}_{\gamma_0}(\mathcal{L}(\Phi))\subseteq\prod_{i=1}^{d}[\lambda_i-a_i, \lambda_{i}+b_{i}]\cap\mathfrak{a}_{+}^{d},
    \end{equation}
        where $a_{i}=2^{i}\gamma_0\delta-\delta$ and $b_{i}=(2^{i}-2)\gamma_0\delta+\delta$.
\end{lemma}
\begin{proof}[Proof of Corollary \ref{CorollaryDeWittGogoglev}]
    Consider $\xi=(\lambda_{1}+\epsilon_{1},\dots,\lambda_{d}+\epsilon_{d} )\not\in \prod_{i=1}^{d}[\lambda_i-a_i, \lambda_{i}+b_{i}]$. Then there is $m\in\{1,\dots, d\}$ such that $\epsilon_{m}<-a_{m}$ or $\epsilon_{m}>b_m$ and $\epsilon_{i}\in [-a_i,b_i]$ for $i<m$.  
    Let $f$ be the linear Schur-convex function defined by $f(\omega)=\sum_{i=1}^{m}\omega_{i}$. 
    
    Then, if $\epsilon_{m}<-a_{m}$
    \begin{align}
    \begin{split}
    f(\xi)&=f(\lambda)+\sum_{i=1}^{m}\epsilon_{i}\\
    &\leq f(\lambda)+\epsilon_{m}+\sum_{i=1}^{m-1}b_{i}\\
    &< f(\lambda)-2^{m}\gamma_0\delta+\delta+\sum_{i=1}^{m-1}((2^{i}-2)\gamma_0\delta+\delta)\\
    &=f(\lambda)-2^{m}\gamma_0\delta+\delta+(2^{m}-2)\gamma_0\delta-\delta(m-1)(2\gamma_0-1)\\
    &=(f(\lambda)+m\delta)-\gamma_0(f(\lambda)+m\delta)+\gamma_0(f(\lambda)-m\delta)\\
    &\leq (1-\gamma_0)\sup_{\mathcal{L}(\Phi)}f+\gamma_0\inf_{\mathcal{L}(\Phi)}f.
    \end{split}
    \end{align}
On the other hand, if $\epsilon_{m}>b_m$
\begin{align}
    \begin{split}
    f(\xi)&=f(\lambda)+\sum_{i=1}^{m}\epsilon_{i}\\
    &\geq f(\lambda)+\epsilon_{m}-\sum_{i=1}^{m-1}a_{i}\\
    &> f(\lambda)+(2^{m}-2)\gamma_0\delta+\delta-\sum_{i=1}^{m-1}(2^{i}\gamma_0\delta-\delta)\\
    &=f(\lambda)+(2^{m}-2)\gamma_0\delta+\delta-(2^{m}-2)\gamma_0\delta+(m-1)\delta\\
    &=f(\lambda)+m\delta\\
    &\geq \sup_{\mathcal{L}(\Phi)}f.
    \end{split}
    \end{align}
In either case, we have shown that $\xi\not\in \mathrm{Aura}_{\gamma_0}(\mathcal{L}(\Phi))$. 
\end{proof}

We can now proceed to prove Corollary \ref{CorollaryDeWittGogoglev}

\begin{proof}[Proof of Corollary \ref{CorollaryDeWittGogoglev}]
    By the previous lemma and by Theorem \ref{MainTheorem} we have that 
    \begin{equation}
    \mathcal{U}(\Phi)\subseteq \mathrm{Aura}_{\gamma_0}(\mathcal{L}(\Phi))\subseteq\prod_{i=1}^{d}[\lambda_i-a_i, \lambda_{i}+b_{i}]\cap\mathfrak{a}_{+}^{d}.
    \end{equation}
    Hence, by Proposition \ref{DominatedSplittingLemma}, to obtain the dominated splitting of index $m$ we have to show that for all $\xi\in  \mathrm{Aura}_{\gamma_0}(\mathcal{L}(\Phi))\subseteq\prod_{i=1}^{d}[\lambda_i-a_i, \lambda_{i}+b_{i}]$, we get that $\xi_{m}>\xi_{m+1}$. 
    
    We see that
    \begin{align}
    \begin{split}
        (2^{m+1}+2^{m}-2)\gamma_0\delta &<\lambda_m-\lambda_{m+1}\\
        \lambda_{m+1}+((2^{m+1}-2)\gamma_0\delta+\delta)&<\lambda_{m}-(2^{m}\gamma_0\delta-\delta)\\
        \lambda_{m+1}+b_{m+1}&<\lambda_{m}-a_{m}.
    \end{split}
    \end{align}
    Hence, we conclude that for for all $\xi\in \mathcal{U}(\Phi)$
    \begin{equation}
    \xi_{m+1}\leq\lambda_{m+1}+b_{m+1}<\lambda_{m}-a_{m}\leq \xi_{m}.
    \end{equation}
    The result follows. 
\end{proof}

Next, we need the following definition to restate and prove \cite[Theorem~1]{Renato}. 
\begin{definition}
    Suppose that $T:X\to X$ is a hyperbolic mixing homeomorphism with constant of expansion $\tau>0$ and $\mathbb{E}$ is a $\theta$-H\"older vector bundle over $X$. We say that a $\theta$-H\"older bundle map $\Phi:\mathbb{E}\to\mathbb{E}$ is \textit{fiber-bunched} if $(\xi_1-\xi_d)^{\top}\leq \theta \tau$.
\end{definition}
Thus, we have the following corollary: 

\begin{corollary} \label{CorollaryVelozo}
    Assume that $T:X\to X$ is a hyperbolic mixing homeomorphism with constant of expansion $\tau>0$, $\mathbb{E}$ is a two dimensional $\theta$-H\"older vector bundle over $X$ and $\Phi:\mathbb{E}\to\mathbb{E}$ is a fiber-bunched  $\theta$-Hölder bundle map covering $T$. Suppose that $\overline{\mathcal{L}(\Phi)}\cap\{\xi\in \mathfrak{a}_{+}^{2}:\xi_1=\xi_2\}=\emptyset$. Then, $\Phi$ has a dominated splitting. 
\end{corollary}

\begin{proof}
    Since $\Phi$ is fiber bunched and $\mathbb{E}$ is a two dimensional vector bundle, by Remark \ref{Aura in dimension 2} and Theorem \ref{MainTheorem}, we have that $\mathcal{U}(\Phi)\subseteq \overline{\mathrm{conv}}(\mathcal{L}(\Phi))$. 
    
    Note that since $\overline{\mathcal{L}(\Phi)}\cap\{\eta:\eta_1=\eta_2\}=\emptyset$, then 
    \begin{equation}
    \min_{\overline{\mathrm{conv}}(\mathcal{L}(\Phi))}(\xi_1-\xi_2)=\inf_{\mathcal{L}(\Phi)}(\xi_1-\xi_2)>0.
    \end{equation}
    Hence, $\mathcal{U}(\Phi)\cap \{\eta:\eta_1=\eta_2\}=\emptyset$. The result follows from Proposition \ref{DominatedSplittingLemma}.
\end{proof}

\section{The uniform spectrum of typical cocycles} \label{Section questions}

The main objective of this paper is to introduce the uniform spectrum as an object worth studying. As it was discussed above, the dominated splittings can be obtained readily from $\mathcal{U}(\Phi)$ and, as we have seen from Theorem \ref{MainTheorem}, there appear to be techniques to bound the set $\mathcal{U}(\Phi)$. The following are a couple of follow up questions that can be posed on the uniform spectrum.    

    We have seen that in general the Lyapunov and uniform spectra can be very different, even in the case of 
    Hölder continuous cocycles over hyperbolic homeomorphisms.  However, it is reasonable to expect that 
    in typical situations these spectra do not differ so much.
    
    In the fiber-bunched setting, pinching and twisting \cite[Definition~1.2]{Avila2006SimplicityOL} are very 
    useful generic properties. Park \cite{Parkproximal,Parkquasimultiplicativity} obtained several results which 
    we can reformulate in terms of the Lyapunov and uniform spectra.  

    \begin{theorem} \label{Park uniform spectrum}
        Assume that $T:X\to X$ is a hyperbolic mixing homeomorphism with constant of expansion $\tau>0$, 
        $\mathbb{E}$ is a $\theta$-H\"older vector bundle over $X$ and 
        $\mathsf{\Lambda}^{j}\Phi:\mathbb{E}\to\mathbb{E}$ is a fiber-bunched, 
        pinching and twisting  $\theta$-Hölder bundle map covering $T$ for all $1\leq j\leq d$.
        Then:
        \begin{enumerate}
            \item $\mathcal{U}(\Phi)=\overline{\mathcal{L}(\Phi)}$.
            \item The uniform spectrum is given by
            \begin{equation} \label{Regularspectrum}
             \mathcal{U}(\Phi)=\left\{\xi\in\mathfrak{a}_{+}^{d}:\text{there exist } x\in X \text{ such that } \xi=\lim_{n\to\infty}\dfrac{1}{n}\kappa^{g}(\Phi^{n}_{x})  \right\}.
            \end{equation}
            \item $\mathcal{U}(\Phi)$ is a convex set. 
        \end{enumerate}
    \end{theorem}
    \begin{proof}
        Under our conditions, Park \cite[Theorem~A]{Parkproximal} states that there exist $C>0$ and $l\in\mathbb{N}$ such that for all $x\in X$ and all $n\in\mathbb{N}$, there exist a periodic point $p\in X$ of of period $n_{p}\in \{n,\dots, n+k\}$ such that
       \begin{equation} \label{TeoremaPark}
        \norm{\kappa^{g}(\Phi^{n}_{x})-\zeta(\Phi^{n_p}_{p})}\leq C.
       \end{equation}

       Since $\mathcal{U}(\Phi)$ is a closed set, it is clear that $\overline{\mathcal{L}(\Phi)}\subseteq \mathcal{U}(\Phi)$. Now, consider $\xi\in \mathcal{U}(\Phi)$. Then there is a sequence $\{x_{j}\}\subseteq X$ and an increasing sequence of natural numbers $\{n_{j}\}\subseteq \mathbb{N}$ such that 
     $$\xi=\lim_{j\to\infty} \dfrac{1}{n_{j}}\kappa^{g}(\Phi^{n_j}_{x_j}).$$
     Using Theorem \ref{TeoremaPark}, we see that for every $x_{j}\in X$ and $n_{j}\in \mathbb{N}$ there is $p_{j}\in X$ of period $m_{j}\in\{n_{j},\dots, n_{j}+k\}$ such that
     \begin{align}
     \begin{split}
         &\lim_{j\to\infty}\norm{\dfrac{1}{n_{j}}\kappa^{g}(\Phi^{n_j}_{x_j})-\dfrac{1}{m_{j}}\zeta(\Phi^{m_j}_{x_j})}\\
         &\hspace{1cm}\leq \lim_{j\to\infty}\norm{\dfrac{1}{n_{j}}\kappa^{g}(\Phi^{n_j}_{x_j})-\dfrac{1}{n_{j}}\zeta(\Phi^{m_j}_{x_j})}+\lim_{j\to\infty}\dfrac{k}{n_{j}}\norm{\dfrac{1}{m_j}\zeta(\Phi^{m_j}_{x_j})}\\
         &\hspace{1cm}\leq \lim_{j\to\infty}\dfrac{1}{n_j}C+\lim_{j\to\infty}\dfrac{k}{n_{j}}\sup_{\eta\in\mathcal{L}(\Phi)}\norm{\eta}\\
         &\hspace{1cm}=0.
        \end{split}
     \end{align}
     Hence, $\xi$ is a limit point of the sequence $\{m_{j}^{-1}\zeta(\Phi^{m_j}_{x_j})\}\subset \mathcal{L}(\Phi)$. Thus, $\xi\in\overline{\mathcal{L}(\Phi)}$. 
     This proves part (1). 

     Let $\mathcal{R}(\Phi)$ be the set in the right hand side of (\ref{Regularspectrum}).
     According to Park \cite[Theorem~D]{Parkproximal}, $\mathcal{R}(\Phi)$ equals the closure of the Lyapunov vectors
     of ergodic measures supported on periodic orbits. The latter set contains $\mathcal{L}(\Phi)$ by \cite[Theorem~1.4]{Kalinin}.
     Therefore, the statement of part (2) follows from part (1). 

    By \cite[Theorem~D]{Parkquasimultiplicativity}, $\mathcal{R}(\Phi)$ is a convex set. Hence part (3) follows.  
     \end{proof}

    Recall that the \textit{relative interior} of a convex set $C$, which we write $\mathrm{rint}(C)$, is the interior 
    of $C$ with respect to its affine hull.  Motivated by \cite[Conjecture~4.7]{MisiurewiczZiemian}, we ask:
    \begin{question}
        Under the hypotheses of Theorem \ref{Park uniform spectrum}, is it true that $\mathrm{rint}(\mathcal{U}(\Phi))\subseteq\mathcal{L}(\Phi)$?
    \end{question}

    If this is the case, we would have that the relative boundary of the uniform spectrum contains all the Lyapunov vectors that do not come from any ergodic measure. 

    We end this paper with the following question about the importance of the fiber bunching condition. 
    \begin{question}
            Are properties (2) and (3) from Theorem \ref{Park uniform spectrum} typical among non-fiber-bunched Hölder bundle maps?
    \end{question}

    \appendix

\section{Some properties of the Lyapunov and uniform spectra} \label{properties Lyapunov and uniform spectra}

\subsection{Proof of Proposition \ref{change of uniform spectrum due to iterations}}

We will need the following elementary ergodic-theoretical lemmas. These results may be known, but due to a lack of precise references we provide the proofs. 

\begin{lemma} \label{existance of eigenvalue}
Let $k\geq 2$ be an integer number, $T:X\to X$ be a continuous map and $\mu\in\mathcal{M}_{\mathrm{erg}}(T)\backslash \mathcal{M}_{\mathrm{erg}}(T^{k})$ be a measure which is ergodic for $T$ but not for $T^{k}$. Then there exists a nonzero function $g\in \mathrm{L}^{2}(\mu)$ and a complex number $\lambda\in\mathbb{C}$ such that $|\lambda|=1$ and $g\circ T=\lambda g$. 
\end{lemma}

\begin{proof}
    Let $U:\text{L}^{2}(\mu)\to \text{L}^{2}(\mu)$ be the Koopman operator of $T$, that is, for any $h\in \text{L}^{2}(\mu)$, $Uh= h\circ T$. Since $\mu\not\in \mathcal{M}_{\mathrm{erg}}(T^{k})$, there exist non-constant $f\in \text{L}^{2}(\mu)$ such that $U^{k}f=f$. Notice that since $\mu\in \mathcal{M}_{\mathrm{erg}}(T)$, $Uf\not= f$. 

    Let $\mathcal{Z}=\mathrm{span}\{f, Uf,\dots, U^{k-1}f\}\subseteq \mathrm{L}^{2}(\mu)$. Note that $\mathrm{Z}$ is a finite dimensional vector space. Note that the restriction $U_{\mathcal{Z}}:\mathcal{Z}\to \mathcal{Z}$ is a unitary operator which is not the identity. Then, by the finite dimensional spectral theorem, $U_{\mathcal{Z}}$ has an eigenvalue $\lambda\not=1.$ Consider a nonzero eigenfunction  $g\in\mathcal{Z}$ such that $Ug=\lambda g$. However, we see that $U_{\mathcal{Z}}=Id$, which implies that $g=U^{k}_{\mathcal{Z}}g=\lambda^{k}g$. The result follows. 
\end{proof}

Now we describe the measures of $\mathcal{M}_{\mathrm{erg}}(T)$ in terms of the measures $\mathcal{M}_{\mathrm{erg}}(T^{k})$.
\begin{lemma} \label{description ergodic measures}
 Suppose that $T:X\to X$ is a continuous map and $k\in\mathbb{N}$ is any positive integer. Then
\begin{equation}
    \mathcal{M}_{\mathrm{erg}}(T)=\left\{\mu=\dfrac{1}{k}\sum_{i=0}^{k-1}T^{i}_{*}\nu \ :\nu\in  \mathcal{M}_{\mathrm{erg}}(T^{k})\right\}.
\end{equation}
\end{lemma}
\begin{proof} First, let us check the $\supseteq$ inclusion. Suppose $\nu\in\mathcal{M}_{\mathrm{erg}}\left(T^{k}\right)$. Since $\nu$ is invariant with respect to $T_{*}^{k}$, we see that
       \begin{align}
           \begin{split}
               T_{*}\left(\dfrac{1}{k}\sum_{i=0}^{k-1}T^{i}_{*}\nu\right)&=\dfrac{1}{k}(T_{*}\nu+\dots+T^{k}_{*}\nu)\\
               &=\dfrac{1}{k}(T_{*}\nu+\dots+\nu)\\
               &=\dfrac{1}{k}\sum_{i=0}^{k-1}T^{i}_{*}\nu,
           \end{split}
       \end{align}
       which shows that $\frac{1}{k}\sum_{i=0}^{k-1}T^{i}_{*}\nu$ is invariant with respect to $T$. 

       Now, consider a set $A\subseteq X$ such that $T^{-1}A=A$. Then, $T^{-m}A=A$ for all $m\in \mathbb{N}$. Since $\nu$ is ergodic with respect to $T^{k}$, we obtain that $\nu(A)=0$ or $\nu(A)=1$. We can conclude that $\frac{1}{k}\sum_{i=0}^{k-1}T^{i}_{*}\nu(A)$ equal to 0 or 1. The inclusion follows. 
    
    Now, we will show the $\subseteq$ inclusion. If $\mu\in\mathcal{M}_{\mathrm{erg}}(T^{k})$, the result follows trivially. 
    
    On the other hand, if $\mu\not\in\mathcal{M}_{\mathrm{erg}}(T^{k})$, consider $|\lambda|=1$ and $g\in \mathrm{L}^{2}(\mu)$ as obtained in Lemma \ref{existance of eigenvalue}. Note that for $\mu$-a.e $x\in X$ $g(Tx)=\lambda g(x)$. Thus, if we define $g^{k}(x)=[g(x)]^{k}$, $g^{k}$ is $T$-invariant and therefore it is $\mu$-a.e constant. This implies that $\mu$-a.e, $g$ attains at most $k$-values. 

    Consider $c\in\mathbb{C}\backslash\{0\}$ such that $\mu(g^{-1}(c))>0$ and denote $A_{0}\coloneq g^{-1}(c)$. Furthermore, for $j\in \{1,\dots, k-1\}$, we define
    \begin{equation}
        A_{j}\coloneq g^{-1}(\lambda^{j}c). 
    \end{equation}
    
    Let $l\in\{1,\dots, k\}$ be the smallest positive integer such that $\lambda^{l}=1$. Then, we can see that $l$ divides $k$, and for $i,j\in \{0,\dots, k-1\}$, either $A_{i}\cap A_{j}=\emptyset$ or $A_{i}=A_{j}$ in the case that $i\equiv j\mod{l}$. 
    
    Note that 
    \begin{equation}\label{invariance of A}
        T^{-1}A_{j}=g^{-1}(\lambda^{j-1}c)=A_{j-1}.
    \end{equation}
    Thus, we see that $\mu\left(\bigcup_{j=0}^{k-1} A_j\right)=1$ since it is a $T$-invariant set which contains $A_{0}$. Hence, for $j\in\{1,\dots, k-1\}$ we can define the measures
    \begin{equation}
        \nu_{j}(B)=\dfrac{\mu(A_{j}\cap B)}{\mu(A_j)},
    \end{equation}
    where $B\subseteq X$ is any $\mu$-measurable set, such that
    \begin{equation}
        \mu(B)=\dfrac{1}{k}\sum_{j=0}^{k-1}\nu_{j}(B).
    \end{equation}

    Additionally, from (\ref{invariance of A}), we see that 
    \begin{equation}
        T_{*}\nu_{j}(B)=\dfrac{\mu(A_{j}\cap T^{-1}B)}{\mu(A_j)}=\dfrac{\mu(T^{-1}(A_{j-1})\cap T^{-1}B)}{\mu(T^{-1}A_{j-1})}=\nu_{j-1}(B),
    \end{equation}
    which implies that 
    \begin{equation}
        \mu=\dfrac{1}{k}\sum_{j=0}^{k-1}T^{j}_{*}\nu_{0}.
    \end{equation}

    It remains to prove that $\nu_{0}\in \mathcal{M}_{\mathrm{erg}}(T^{k})$. First, since $T^{-k}A_{0}=A_{0}$, we obtain that $\nu_{0}$ is $T^{k}_{*}$-invariant. Now, if $T^{-k}C=C$ $\nu_{0}$-a.e, without loss of generality we can assume that $C\subseteq A_{0}$. Hence, $D=\bigsqcup_{j=0}^{l-1}T^{-j}C$ is a $T$-invariant set. Thus, $\mu(D)=0$ or $\mu(D)=1$, which implies that $\nu_{0}(C)=0\text{ or }1$. The result follows.   \qedhere
\end{proof}

We can now show the main proof of this subsection. 

\begin{proof}[Proof of Proposition \ref{change of uniform spectrum due to iterations}] We begin by showing the equality (\ref{Lyapunov spectrum iterations}). Define, for each $x\in X$,
\begin{equation}
    \lambda(x)\coloneq \limsup_{m\to \infty}\dfrac{1}{m}\kappa^{g}(\Phi^{m}_{x}),
\end{equation}
    where the limsup is taken entrywise. Note that the Borel function $\lambda$ is $T$-invariant: this can be shown using the singular value bounds \cite[(4.6)]{BPS26}. By the Furstenberg-Kesten theorem, 
    \begin{equation}
        \mathcal{L}(\Phi)=\left\{\int\lambda d\mu: \ \mu\in\mathcal{M}_{\mathrm{erg}}(T)\right\}.
    \end{equation}
    
    For every $k>0$, 
    \begin{equation}
        \limsup_{n\to\infty} \dfrac{1}{nk}\kappa^{g}(\Phi^{kn}_{x})=\lambda(x),
    \end{equation}
    again as a consequence of the bounds \cite[(4.6)]{BPS26}. In particular, $\frac{1}{k}\mathcal{L}(\Phi^{k})$ is the set of integrals $\int \lambda d\nu$, where $\nu$ runs on the set $\mathcal{M}_{\mathrm{erg}}(T^{k})$. By the Lemma \ref{description ergodic measures}, for every $\nu\in \mathcal{M}_{\mathrm{erg}}(T^{k})$, the measure $\mu=\frac{1}{k}\sum_{i=0}^{k-1}T^{i}_{*}\nu$ is ergodic for $T$ and, conversely, every $\mu\in \mathcal{M}_{\mathrm{erg}}(T)$ can be written in the same way as above for some $\nu\in \mathcal{M}_{\mathrm{erg}}(T^{k})$. Note that since $\lambda$ is $T$-invariant, then $\int \lambda d\mu=\int \lambda d\nu$. It follows that $\mathcal{L}(\Phi)=\frac{1}{k}\mathcal{L}(\Phi^{k})$, as we expected. 

     Now, we will focus on the equality (\ref{Uniform spectrum iterations}). Suppose that $\xi\in \mathcal{U}(\Phi^{k})$. Then there are sequences $\{x_j\}\subseteq X$ and $\{n_j\}\subseteq \mathbb{N}$ where $n_j\to\infty$ such that
    \begin{equation}
        \xi=\lim_{j\to\infty}\dfrac{1}{n_j}\kappa^{g}(\Phi^{kn_j}_{x_j})=\lim_{j\to\infty}\dfrac{k}{kn_j}\kappa^{g}(\Phi^{kn_j}_{x_j})=k \omega,
    \end{equation}
    where $\omega\in \mathcal{U}(\Phi)$. Hence $\mathcal{U}(\Phi^{k})\subseteq k~\mathcal{U}(\Phi).$

    Now suppose that $\omega\in \mathcal{U}(\Phi)$. Then there are sequences $\{y_j\}\subseteq X$ and $\{m_j\}\subseteq \mathbb{N}$ where $m_j\to\infty$ such that $\omega=\lim_{j\to\infty}\frac{1}{m_j}\kappa^{g}(\Phi^{m_j}_{y_j}).$

    Take $\delta>0$. Using the bounds \cite[(4.6)]{BPS26}, for $m\in\mathbb{N}$ sufficiently large and all $x\in X$
    \begin{equation} \label{Convergence of Cartan projection}
        \left\| \frac{1}{m} \kappa^{g}(\Phi^m_x) -  \frac{1}{m+1} \kappa^{g}(\Phi^{m+1}_{x}) \right\|<\delta.
    \end{equation}
    Thus, for each $m_{j}$ large enough, we can find $l_j\in\mathbb{N}$ such that $m_{j}\leq kl_j\leq m_j +k$ so that 
    \begin{align}
    \begin{split}
        \norm{\omega-\frac{1}{ kl_j}\kappa^{g}(\Phi^{k l_j}_{y_j})}&\leq \norm{\omega-\frac{1}{m_j}\kappa^{g}(\Phi^{m_j}_{y_j})}+\norm{\frac{1}{m_j}\kappa^{g}(\Phi^{m_j}_{y_j})-\frac{1}{ kl_j}\kappa^{g}(\Phi^{k l_j}_{y_j})}\\
        &< (k+1) \delta. 
        \end{split}
    \end{align}
    Hence, we conclude that $\omega \in \frac{1}{k}\mathcal{U}(\Phi^{k})$. Hence, we obtain the equality $\ref{Uniform spectrum iterations})$. 

    The equalities (\ref{Inverse Lyapunov Spectrum}) and (\ref{inverse uniform spectrum}) follow from the following observations: if $T:X\to X$ is invertible, then $\mathcal{M}_{\mathrm{erg}}(T)=\mathcal{M}_{\mathrm{erg}}(T^{-1})$ and for every $A\in\mathrm{GL}(\mathbb{E}_{x},\mathbb{E}_{y})$ where $x,y\in X$ then $\kappa^{g}(A^{-1})=\iota(\kappa^{g}(A))$.  
\end{proof}

\subsection{Proof of Proposition \ref{p.criterion_U_connected}}

The proof of Proposition~\ref{p.criterion_U_connected} will need the following lemma:

\begin{lemma}\label{l.connect}
If $T$ is invariantly connected, then for any $\varepsilon>0$ and any two points $x,y \in X$, we can find a list of points $x_0,x_1,\dots,x_k$ in $X$ with $x_0=x$, $x_k=y$, and for every $i \in \{0,\dots,k-1\}$, either $d(x_i,x_{i+1}) < \varepsilon$, or $x_{i+1} = T(x_i)$, or $x_i = T(x_{i+1})$.
\end{lemma}

\begin{proof}
Given $\varepsilon>0$ and $x \in X$, let $U_0 \coloneq \{x\}$ and define recursively $U_{i+1} \coloneq N_\varepsilon(U_i) \cup T^{-1}(U_i) \cup T(U_i)$, where $N_\varepsilon$ denotes the $\varepsilon$-neighborhood. 
Let $U \coloneq \bigcup_{i=0}^\infty U_i$.
Then $N_\varepsilon(U) = U$; in particular, $U$ is a clopen set. 
Furthermore, $T^{-1}(U) = U$.
Since $T$ is invariantly connected, $U=X$. 
The statement of the lemma immediately follows.
\end{proof}

\begin{proof}[Proof of Proposition~\ref{p.criterion_U_connected}]
Assume for a contradiction that the set $\mathcal{U}(\Phi)$ is not connected. 
Then $\mathcal{U}(\Phi) \subseteq V \cup W$ for a pair of disjoint open sets $V, W \subseteq \mathbb{R}^d$, both of which have nonempty intersection with $\mathcal{U}(\Phi)$.
Fix $\delta>0$ such that $\|\xi - \eta\|>\delta$ for all $\xi \in V\cap \mathcal{U}(\Phi)$ and $\eta \in W\cap \mathcal{U}(\Phi)$.
Since $\mathcal{U}(\Phi) = \mathsf{Ls}_{n \to \infty} \frac{1}{n} \kappa^{g}(\Phi^n)$, we have
$\frac{1}{n} \kappa^{g}(\Phi^n) \subseteq V \cup W$ for all sufficiently large $n$.

On the other hand, by \cite[(4.6)]{BPS26}, for $n$ large enough and for all $x\in X$
\begin{equation}\label{e.slow_variations}
\left\| \frac{1}{n} \kappa^{g}(\Phi^n_x) -  \frac{1}{n} \kappa^{g}(\Phi^n_{Tx}) \right\| < \delta \quad.
\end{equation}

It follows that for every sufficiently large $n$, the set $\frac{1}{n} \kappa^{g}(\Phi^n)$ is contained either in $V$ or in $W$. To verify this, consider $x,y\in X$ and apply Lemma~\ref{l.connect} with $\varepsilon = \varepsilon_n$ small enough. Using (\ref{e.slow_variations}), there are $x_{0},x_1,\dots, x_k$ where $x_0=x$, $x_k=y$ such that  
\begin{equation}
\left\| \frac{1}{n} \kappa^{g}(\Phi^n_{x_i}) -  \frac{1}{n} \kappa^{g}(\Phi^n_{x_{i+1}}) \right\| < \delta. 
\end{equation}
Hence, both $\frac{1}{n} \kappa^{g}(\Phi^n_{x}), \frac{1}{n} \kappa^{g}(\Phi^n_{y})$ belong either to $V$ or to $W$. 

On the other hand, since $\mathcal{U}(\Phi)$ intersects $V$ (resp.~$W$), the sequence of sets $\frac{1}{n} \kappa^{g}(\Phi^n)$ must intersect $V$ (resp.~$W$) infinitely often.
It follows that there exists $n \gg 1$ such that $\frac{1}{n} \kappa^{g}(\Phi^n) \subseteq V$ and $\frac{1}{n+1} \kappa^{g}(\Phi^{n+1}) \subseteq W$.
Then the inequality \eqref{Convergence of Cartan projection} is violated. 
This contradiction finishes the proof.
\end{proof}

\section{Exponential closing lemma} \label{Appendix B}


Let us develop the tools to prove the Lemma \ref{approximation}. First, for $m\geq 1$, we consider the following subsets of the global stable and unstable set respectively
\begin{align}
\begin{split}
    W^{s}_{m}(x)&\coloneq T^{-m}W^{s}_{\epsilon_{0}}(T^mx),\\
    W^{u}_{m}(x)&\coloneq T^{m}W^{u}_{\epsilon_{0}}(T^{-m}x). 
\end{split}
\end{align}

Now we proceed to show that iterations of the local stable and unstable sets are $\delta$-close to any point on the metric space. 
\begin{lemma} \label{densidadEstable}
    Suppose $(X,d)$ is a compact metric space and $T:X\to X$ is a topologically mixing hyperbolic homeomorphism. Then for all $x \in X$ and $0<\delta\leq \epsilon_{1}$ there exist $N>0$ such that $T^{-n}W^{s}_{\epsilon_{1}}(T^{n}x)$ and $T^{n}W^{u}_{\epsilon_{1}}(T^{-n}x)$ are $\delta-$dense in $X$ for all $n\geq N$.
\end{lemma}
\begin{proof}
    Consider $x\in X$. Since $X$ is compact, there exist $z_{1},\dots, z_{m}\in X$ such that $X\subseteq \bigcup_{i=1}^{m}B(z_{i},\delta/4C)$ where $C>1$ is as in the definition of hyperbolic homeomorphism. Due to topological mixing, we can find $N(\delta)\in\mathbb{N}$ such that $T^{n}B(z_{i},\delta/4C)\cap B(z_{j}, \delta/4C)\not=\emptyset$ for $1\leq i, j\leq m$ and all $n\geq N$. We see that there exist some $z_{k}$ such that $T^{n}(x)\in B(z_{k}, \delta/4C)$. 
    
    For every $z_{i}$ there is $w\in B(z_{i}, \delta/4C)$ such that $T^{n}(w)\in B(z_{k}, \delta/4C)$. Denote $T^{n}(w')=[T^{n}(w),T^{n}(x)]$. Thus $T^{n}(w')\in W^{u}_{\epsilon_{0}}(T^{n}(w))$, so we see that $d(w',w)<\delta/2$. Hence we have that $d(z_{i},w')\leq \delta$.

    We conclude by noting that $w'\in T^{-n}W^{s}_{\epsilon_{1}}(T^nx)$. A similar proof yields the same result for $T^{n}W^{u}_{\epsilon_{1}}(T^{-n}x)$.
\end{proof}
Using the density of the stable and unstable sets, we can show that their intersection is nonempty.
\begin{lemma} 
    Under the hypotheses of Lemma \ref{densidadEstable}, there exist $N>0$ so that for every $x,y\in X$ we have that $W^{s}_{N}(x)\cap W^{u}_{N}(y)\not=\emptyset$.
\end{lemma}
\begin{proof}
    Consider $N_{1}>0$ such that $e^{-N_{1}\tau}<\epsilon_{0}-\epsilon_{1}$.  From the previous lemma, we know that there is $N_{2}>0$ such that $T^{-n}W^{s}_{\epsilon_{1}}(T^{n}x), T^{n}W^{u}_{\epsilon_{1}}(T^{-n}x)$ are $\epsilon_{1}/2$-dense for all $n\geq N_{2}$. Take $N=\max\{N_1,N_2\}$. Hence, there exist $T^{-N}y'\in W^{u}_{\epsilon_1}(T^{-N}y)$, $T^{N}x'\in W^{s}_{\epsilon_{1}}(T^Nx)$ such that $y'$  and $x'$ are sufficiently close so as to define $z=[y',x']$. We see that $z\in W^{s}_{\epsilon_{0}}(x')$. Thus, we have the following bounds for any $n\geq N$
    \begin{align}
    \begin{split}
        d(T^{n}z,T^{n}x)&\leq d(T^{n}z,T^{n}x')+d(T^{n}x',T^{n}x)\\
        &\leq e^{-\tau n}d(z,x')+e^{-\tau(n-N)}d(T^{N}x',T^{N}x)\\
        &\leq e^{-\tau N}+\epsilon_{1}\\
        &\leq \epsilon_{0}.
    \end{split}
    \end{align}
    We conclude that $z\in W^{s}_{N}(x)$. Via a similar procedure, we can find that $z\in W^{u}_{N}(y)$. Hence, the intersection is nonempty. 
\end{proof}

The next lemma ensures that every orbit segment can be approximated by a periodic pseudo-orbit for $n$ sufficiently large. 
\begin{lemma}  \label{pseudoorbit}
    Under the conditions of Lemma \ref{densidadEstable}, for all $x\in X$ and $\epsilon>0$, there is $R>0$ and $n_{0}(\epsilon)\in \mathbb{N}$ such that for $n\geq n_{0}$ there is a periodic $\epsilon-$pseudo-orbit $\{z_0,\dots, z_{n}\}$ such that for $0\leq i\leq n$
    \begin{align}
      d(T^{i}x, z_i)\leq Re^{-\tau \min\{i,n-i\}}.
    \end{align}
\end{lemma}
\begin{proof}
    Consider $\epsilon>0$. From the previous lemma, we can fix $N>0$ such that there exist $z\in W^{s}_{N}(x)\cap W^{u}_{N}(T^{n}x)$ for $n\gg 1$. Since $T^{N}z\in W^{s}_{\epsilon_{0}}(T^{N}x)$, we see that for $N\leq m\leq n/2$ 
    \begin{align}
        d(T^{m}x, T^{m}z)&\leq e^{-\tau(m-N)}d(T^{N}x,T^{N}z),\\
        \intertext{and since $T$ is Lipschitz and $N$ is fixed}
        &\leq \text{Lip}(T)^{N} e^{\tau N}e^{-\tau m}d(x,z).
        \intertext{Similarly, since $T^{-N}z\in W^{u}_{\epsilon_{0}}(T^{n-N}x)$}
        \begin{split}
        d(T^{n-m}x, T^{-m}z)&\leq e^{-\tau(m-N)}d(T^{n-N}x,T^{-N}z)\\
        &\leq \text{Lip}(T^{-1})^{N} e^{\tau N}e^{-\tau m}d(z,T^nx).
        \end{split}
    \end{align}
    Hence, it is possible to find $R>0$ such that for $0\leq m\leq \floor{n/2}$
    \begin{align}
    \begin{split}
      d(T^{m}x, T^{m}z)&\leq Re^{-\tau m} \\
      d(T^{n-m}x, T^{-m}z)&\leq Re^{-\tau m}. 
      \end{split}
    \end{align}
    
    Since $T^{\floor{n/}2}x=T^{n-\ceil{n/2}}x$, we get that for $n_0\in \mathbb{N}$ sufficiently large and $n\geq n_0$
    \begin{align}
    \begin{split}
      d(T^{\floor{n/2}}z, T^{-\ceil{n/2})}z)&<Re^{-\tau \floor{n/2}}+Re^{-\tau \ceil{n/2}}\\
      &<\epsilon.
    \end{split}
    \end{align}
    Thus, if we allow $z_{i}=T^{i}z$ for $0\leq i\leq \floor{n/2}-1$ and $z_{i}=T^{-(n-i)}z$ for $\floor{n/2}\leq i\leq n$ we obtain a $\epsilon$-pseudo-orbit which yields the expected results.  
\end{proof}

Now we have all the tools to prove Lemma \ref{approximation}.

\begin{proof}[Proof of Lemma \ref{approximation}]
    By the Lipschitz shadowing lemma for hyperbolic homeomorphisms \cite[Theorem~2]{Sakai}, there are $\delta>0$ and $L>0$ such that for $0<\epsilon<\delta$ and any $\epsilon$-pseudo-orbit $\{z_{0},\dots, z_{n}\}$ there is $p\in X$ such that $d(z_{i},T^{i}y)<L\epsilon$ and $p=T^{n}p$. 

    By Lemma \ref{pseudoorbit}, we can find $R>0$ and $n_{0}\in\mathbb{N}$ such that for all $n\geq n_{0}$ and all $x\in X$ there is a periodic $Re^{-\tau \floor{n/2}}$ pseudo-orbit $\{z_{0},\dots, z_{n}\}$ which approximates the orbit segment $x,\dots, T^{n}x$. Note that by enlarging $n_0$ we may consider $\epsilon=Re^{-\tau \floor{n/2}}<\delta$.

    Consequently, by the Lipschitz shadowing lemma, we obtain $p\in X$ which approximates $\{z_{0},\dots, z_{n}\}$ and has the previously discussed properties. 
    
    The result is obtained by observing that there is $C>0$ such that
    \begin{align}
    \begin{split}
        d(T^{i}x, T^{i}p)&\leq d(T^{i}x,z_{i})+d(z_{i},T^{i}p)\\
        &\leq Ce^{-\tau\min\{i,n-i\}}.
    \end{split}
    \end{align}
    We remark that the constants obtained here do not depend on the chosen point $x$ or on the iterations $n$. 
    \end{proof}

\bibliographystyle{alpha}
\bibliography{LyapunovAndUniformSpectra}
\end{document}